\documentclass[final,11pt,3p]{article}%
\usepackage{setspace}
\usepackage{amsfonts}
\usepackage{amssymb}
\usepackage{amsthm}
\usepackage{amsmath}
\usepackage{algorithmic}
\usepackage{color}
\usepackage{graphicx}
\usepackage{verbatim}
\usepackage{subfigure}
\usepackage[toc,page]{appendix}%
\providecommand{\U}[1]{\protect\rule{.1in}{.1in}}

\newtheorem{theorem}{Theorem}

\newtheorem{definition}{Definition}

{\theoremstyle{THkey}

}

\newcommand{\citealp}{\cite}
\newcommand{\halmos}{}
\newcommand{\Halmos}{}
\providecommand{\keywords}[1]{{\textit{Key words:}} #1}

\newcommand{\review}[1]{{\color{black}#1}}

\begin{document}


\title{Rare-Event Simulation for Distribution Networks}
\author{Jose Blanchet \\Department of Industrial Engineering and Operations Research, Department of Statistics\\ Columbia University, New York, NY 10027, jose.blanchet@columbia.edu\\\vspace{0.2cm}
\and Juan Li \\Department of Industrial Engineering and Operations Research\\Columbia University, New York, NY 10027, jl3035@columbia.edu\\\vspace{0.2cm}
\and Marvin K. Nakayama\\Department of Computer Science, New Jersey Institute of Technology\\ University Heights Newwark, NJ 07102, marvin@njit.edu
}
\maketitle

\begin{abstract}
\noindent We model equilibrium allocations in a distribution network as the solution of
a linear program (LP) which minimizes the cost of unserved demands across nodes in the network. The constraints in the LP dictate that once a given
node's supply is exhausted, its unserved demand is distributed among
neighboring nodes. All nodes do the same and the resulting solution is the
equilibrium allocation. Assuming that the demands are random (following a
jointly Gaussian law), our goal is to study the probability that the optimal
cost (i.e. the cost of unserved demands in equilibrium) exceeds a large
threshold, which is a rare event. Our contribution is the development of
importance sampling and conditional Monte Carlo algorithms for estimating this
probability. We establish the asymptotic efficiency of our algorithms and also
present numerical results which illustrate strong performance of our procedures.
\end{abstract}


\noindent\keywords{ distribution network; linear program; rare event simulation; importance sampling; conditional Monte Carlo }


%


\section{Introduction}

Consider the following model of a distribution network. We assume that there
is a commodity to be distributed among various nodes in a network. Each node
is endowed with a given supply of the commodity and at the same time it
experiences a random demand. We assume that the commodity is infinitely
divisible. If the demand at a given node exceeds its supply, then the excess
demand is distributed according to some proportions to each of its neighbors,
which in turn do the same. In order to obtain the distribution amounts in
equilibrium, we solve a linear program (LP), where the objective function to
minimize is the sum across nodes of the unserved demands.

One possible practical example where such a problem might arise is an
electric power grid. Each node represents a geographic region, and there is an edge between two nodes if transmission lines directly connect them. Each region has generators, which provide the
region's supply of electricity. Also, each region has a random load (i.e.,
demand for electricity) from consumers. If a
region's load exceeds its supply, then the network tries to serve a node's
excess load by sending it to neighboring regions. One of the most important issues in operating a power grid is to keep the stability of the network and make sure demands can be satisfied. 
If the total amount of load
not served at their originating regions exceeds a threshold $k$, then we
consider the network to have failed. To better operate this power transmission system, it is essential to estimate the probability that this network fails.

Another application involves load distribution for internet services, such as
web servers, cloud-computing services, and domain name servers (DNS). A
company may have a number of fixed-capacity servers situated in different
geographic regions. As the requests to servers (i.e. the demand) arrive, a
specific server tries to fulfill its own local requests, but if the demand
exceeds its capacity, then the server may offload its excess to a neighboring
server. Since this shifting may incur additional delays for the user, we want
to minimize the amount of distributed demand. This is similar to load
balancing; e.g., see \cite{k}.

Let $\alpha\left(  k\right)  $ be the probability that the sum of unserved
demands, in equilibrium, exceeds threshold $k$. \review{Our goal is to estimate the
probability $\alpha\left(  k\right)  $, with $k=k_{n}$, where $n$ is a rarity
parameter and we scale the supply as a function of $n$ and we let $n$ increase
reflecting a situation in which the supply is large. The parameter $k=k_{n}$
is allowed to grow with $n$ or can be constant. }Assuming
jointly distributed multivariate Gaussian demands, we provide asymptotically
optimal estimators, together with numerical experiments showing their
performance, and associated large deviations results. We recall that an
unbiased estimator for $\alpha\left(  k_n\right)  $ is \textit{asymptotically optimal} when $n$ goes to infinity if
the logarithm of its second moment is asymptotically equivalent to the
logarithm of $\alpha ^2\left(  k_n\right) $ (see \citealp{ag}, for
notions of efficiency in rare-event simulation).

As far as we know, this paper provides the first type of large deviations
analysis and efficient Monte Carlo for solutions of linear programs with
random input. More precisely, our contributions are as follows:

1) For our model formulation, we show that our optimal allocation
is invariant if one replaces the objective function by any other criterion that is increasing as a function of the unserved demands (see Theorem \ref{thm_insensitivity}).

2) We establish large deviations analysis for our class of linear
programs with random input (see Theorem \ref{same_order}).

3) We develop an importance sampling (IS) algorithm for estimating
$\alpha\left(  k_n\right)  $, and we show that the algorithm is asymptotically
optimal as the supply gets large, \review{and the threshold $k_n$ is a constant or increases with $n$} (see Section~\ref{sec:IS}).

4) We develop a conditional\ Monte Carlo (CMC) algorithm for the evaluation of
$\alpha\left(  k_n\right)  $, and we prove the asymptotic optimality of this
procedure as the supply gets large, \review{and the threshold $k_n$ is a constant or increases with $n$} (see Section~\ref{sec:CMC}).

5)\ We provide several numerical examples in Section~\ref{sec:exp} that validate the performance of
our algorithm.

Some of the results regarding CMC previously appeared in a conference version of this paper
(\citealp{bln}). Our conference paper restricted the LP's objective function to be the sum of the unserved demands, and we now prove its invariance, as described in contribution 1), which greatly expands the applicability of our approach. Regarding contribution 2), we study the asymptotic behaviors of this network which is not discussed in the conference version. Regarding contribution 3), we develop an importance sampling algorithm which is not studied in the conference version, and we provide a proof of asymptotic optimality and algorithm implementation. As for contribution 4), although in the conference version, we have studied the CMC algorithm and its implementation (see \citealp{bln}, Section~4.3), no mathematical proof is provided regarding the asymptotic optimality of this algorithm. Here, in the journal version, we prove it rigorously. Finally, regarding contribution 5), instead of only comparing the naive simulation and CMC, we compare IS as well. In addition, to show the asymptotic optimality of our algorithms, we include numerical examples in which the rarity parameter changes.

We now explain how our paper relates to prior work. First, regarding 1), we note that similar
results, with different types of networks and other applications, have been
obtained in the literature (see \citealp{en}). We only learned about these applications after we obtained our
model formulation, but we believe the connections are relevant. For the IS algorithm (contribution 2), we introduce a probability
measure that is obtained by connecting the event of interest (i.e. total
unserved demands in equilibrium exceeding a threshold) with a simple union
event involving the demands. Then we use an IS distribution inspired by an approach developed by \cite{abl}. \review{IS algorithms have also been used in \cite{wad}, and \cite{pls} to solve a network operation problem with random inputs. While those two papers focus on the assessment of electrical constraints violation, }we make use of IS to assess the solution of LP which involves optimization. Regarding the CMC estimator, we express
the Gaussian demands in polar coordinates. Given the angle, the conditional
probability of the LP's optimal objective function value exceeding $k$ can be
expressed as the probability of the radial component of the Gaussian lying in
an interval or union of intervals, and this conditional probability can be
computed analytically. The use of polar transformations for CMC and rare event simulation has been
used in the past, see for example, \cite{a}. \cite{ag}, Chapters V
and VI, provide additional background material on importance sampling and conditional
Monte Carlo.

Our work also has other potential
applications, in particular to cascading failures, which has been an interesting and important research topic. For example,
\cite{w} studies cascades in a sparse, random network of interacting
agents whose decisions are determined by the actions of their neighbors
according to a simple threshold rule. \cite{d} consider a branching
process model of cascading failures in an electric power grid. \cite{ing}
analyze a continuous-time Markov chain of a dependability model with cascading failures.
\review{ For another example, \cite{bbl} study the temperature evolution of transmission lines. While that paper provides control algorithms to limit the probability that
cascading failures happen. Our IS strategy (discussed in Section \ref{sec:IS} ) can be
modified to estimate the probability of observing a cascading failure under the
policy advocated by that paper. The modification consists of defining, for
instance, the objective function to be minimized as the worst case
temperature over the lines in the network. Additional constraints need to be
modified or approximated by linear constraints. \cite{wetal} also study rare-event simulation for analyzing blackouts.}

We would like to point out that although we assume multivariate Gaussian demands in this paper, the CMC algorithm can be applied to the case when the demands follow an elliptical distribution (see \citealp{mfe}). Furthermore, while an elliptical copula exhibits symmetric tail dependence, the well known Archimedean copula allows asymmetric tail dependence (\citealp{bhc}). Making use of the results in \cite{mn}, we can see that CMC algorithm is also applicable to Archimedean copula, which makes this algorithm very powerful in solving a wide range of problems. also study rare-event simulation for analyzing blackouts.

The rest of the paper develops as follows.
Section~\ref{sec:model} presents the model of the distribution network, and it
also defines the LP problem and its dual. We establish some properties of the
primal and dual LPs in Section~\ref{sec:properties}. The asymptotic behavior
of the model is discussed in Section~\ref{sec:asym}. We describe the
asymptotic optimality and implementations of importance sampling and
conditional Monte Carlo methods for estimating $\alpha(k_n)$ in
Section~\ref{sec:sim}. Section~\ref{sec:exp} contains the experimental results
from some examples, and we give some final comments in
Section~\ref{sec:conc}.

\section{Model Description}

\label{sec:model}


As we introduce our model and discuss its properties we will follow closely
the discussion in \cite{bln}. Suppose there is a directed graph $G=(V,E)$, where $V=\{1,2,\dots,d\}$ is the set of
vertices and $E=\{(i,j):\exists~\text{directed edge from vertex}%
~i~\text{to}~\text{vertex}~j\}$ is the set of edges. The incidence matrix of
the graph is denoted by $H=(H(i,j):i,j\in V)$, where $H(i,j)=1$ if $(i,j)\in
E$, and $H(i,j)=0$ otherwise, and we assume $H(i,i)=0$ for any $i\in V$. The network model we consider is induced by this graph, and we also assume the following:
\begin{description}
\item {\textbf{1}} The network is irreducible in the sense that the matrix $H$ is irreducible.

\item {\textbf{2}} Each node $i$ has a given fixed supply $s_{i}$.

\item {\textbf{3}} Each node $i$ is subjected to a random demand $D_{i}$. The demand vector $\boldsymbol{D} = (D_{1}, D_{2}\dots
,D_{d})^{\prime}$ is
jointly Gaussian $N(\boldsymbol{\mu}, \Sigma)$, where
prime denotes transpose, $\boldsymbol{\mu}$ is the mean vector, and $\Sigma$ is the covariance matrix. 

\item {\textbf{4}} The expectation of $D_{i}$ is less than or equal
to $s_{i}$ for each node $i$.
\end{description}

Each node tries to serve its realized demand. However, if a given node's supply is exhausted, it distributes the unserved demand to its neighbors,
which, in turn, do the same with their respective neighbors. Nevertheless, there is a cost associated with transferring unserved demands which should be minimized. We construct a linear program to describe this problem. The demands achieve an equilibrium point at each feasible solution, and the objective function is to minimize the sum of the excess demands across the nodes. Let $\boldsymbol{s} = (s_{1},s_{2},\dots,s_{d})^{\prime}$, and the LP is:
\begin{align}
\label{original_primal}
\min~~  &  \sum_{i=1}^{d} x_{i}^{+}\nonumber\\
\text{s.t.}~~  &  D_{i} - s_{i} + \sum_{j:(j,i)\in E}x_{j}^{+} a_{ji}%
=x_{i}^{+}-x_{i}^{-} , \forall i\nonumber\\
&  x_{i}^{+} \ge0, x_{i}^{-} \ge0, \forall i.
\end{align}

The quantity $x_{i}^{+} \ge0$ represents the shedded demand from node $i$ in
equilibrium, which is distributed among its neighbors using a fixed
distribution scheme, which we describe shortly. The quantity $x_{i}^{-} \ge0$
represents the unused supply at node $i$ in equilibrium. Therefore, in
equilibrium, if $x_{i}^{+} - x_{i}^{-} > 0$, then node $i$ sheds demand; if
$x_{i}^{+} - x_{i}^{-} < 0$, then node $i$ has unused supply. When node $j$
has excess demand, $a_{ji}$ denotes the proportion of unserved demand at node
$j$ distributed to node $i$. We assume that if $H(i,j) = 0$, then $a_{ij} =
0$; if $H(i,j) = 1$, then $a_{ij} > 0$. In addition, $\sum_{j = 1}^{d} a_{ij}
= 1, \forall i = 1,2,\dots, d$. The solution moves around excess demands and
supplies to neighbors but does so in such a way that the sum of $x_{i}^{+}$'s,
which are the equilibrium shedded demands, is minimized. The problem can be
expressed in matrix notation as follows. Define $A(i,j)=a_{ij}$ (note that
$A(i,i) = 0$). Let $\boldsymbol{1} = (1,1,\dots,1)^{\prime}$ denote the
$d$-dimensional column vector with all components equal to $1$. Then the
previous linear programming problem (\ref{original_primal}) can be written as:
\begin{align}
\label{primal}\min~~  & \boldsymbol{1}^{\prime}\boldsymbol{x}^{+} +
\boldsymbol{0}^{\prime}\boldsymbol{x}^{-}\nonumber\\
\text{s.t.}~~ &  (A^{\prime}- I) \boldsymbol{x}^{+} + I \boldsymbol{x}^{-} =
\boldsymbol{s}-\boldsymbol{D}\nonumber\\
&  \boldsymbol{x}^{+} \ge\boldsymbol{0}, \boldsymbol{x}^{-} \ge\boldsymbol{0},
\end{align}
where $\boldsymbol{0} = (0,0,\dots,0)^{\prime}$ is the $d$-dimensional column
vector with all components equal to $0$, $A = (A(i,j): i,j \in V)$, $I$ is the
$d \times d$ identity matrix, $\boldsymbol{x}^{+} = (x_{1}^{+},x_{2}^{+}%
,\dots,x_{d}^{+})^{\prime}$, and $\boldsymbol{x}^{-} = (x_{1}^{-},x_{2}%
^{-},\dots,x_{d}^{-})^{\prime}$. The goal is that the sum of shedded demands
is as small as possible because, e.g., the cost of distributing demands is
high. If the cost is too high, for example, larger than a given number, say
$k$, or the LP is infeasible, we consider the network to have failed.

Note that while in \cite{bln}, we assume that the unserved demands are equally distributed to neighbors, here we make a small but important extension. We allow the proportions to be any non-negative numbers.

Now, we also introduce the dual linear program:
\begin{align}
\label{dual}\max~~  & \boldsymbol{y}^{\prime}\boldsymbol{r}\nonumber\\
\text{s.t.}~~  &  M\boldsymbol{y} \le\boldsymbol{1}\nonumber\\
&  \boldsymbol{y} \ge\boldsymbol{0},
\end{align}
where $M = I - A$ and $r = \boldsymbol{D}-\boldsymbol{s}$.

We are interested in computing the probability that the network fails, for different
values of $k$. Let $\alpha(k)$ represent this failure probability, and $L(\boldsymbol{D})$ denote the optimal value of the dual when the demand
vector is $\boldsymbol{D}$. As discussed in \cite{bln},
\begin{align}
\label{alpha}\alpha(k)  = \beta_{0} + \beta_{1}(k) = P\{L(\boldsymbol{D}) > k \},
\end{align}
where $\beta_{0}$ is the
probability that the primal is infeasible, and $\beta_{1}(k)$ is the probability that the primal is feasible, but the cost is larger than $k$. 

\review{Since the discussion in Section \ref{sec:properties} is valid for all $k$, we do not define $k$ as a function of the rarity parameter $n$ until Section \ref{sec:asym} }.

\section{ Properties of Our Primal and Dual Linear Programs}

\label{sec:properties}


\subsection{Feasibility of the Solutions to the Primal and Dual}
Our previous conference paper proves two theorems on properties of the primal and dual LPs for the special case when $A(i,j) = H(i,j)/\sum_{l=1}^d H(i.l)$. We claim that both theorems are still valid for our more general $A(i,j)$, and the proofs are exactly the same. Here we only list the property regarding feasibility which will be used later, but omit the proof.

\begin{theorem} {\color{white}1}\label{thm_feasibility} 
\begin{description}

\item[(a)] The dual problem (\ref{dual}) is always feasible.

\item[(b)] The primal problem (\ref{primal}) is feasible if and only if
$\sum_{i=1}^{d} D_{i}\le\sum_{i=1}^{d} s_{i}$.
\end{description}
\end{theorem}

\subsection{Uniqueness and Positivity of the Solution to the Primal}

\begin{theorem}
\label{thm_positive} When the primal problem (\ref{primal}) is feasible, it has the following properties:

\begin{description}

\item[(a)] It has a unique optimal solution.

\item[(b)] At the optimal solution, at most one element in the pair
$(x_{k}^{+}, x_{k}^{-})$ is strictly positive, $\forall1 \le k \le d$.
\end{description}
\end{theorem}

To emphasize the main results of the paper, we postpone the formal proof to Appendix \ref{appendix: thm_positive}, and only give a brief explanation here. For (a), assuming there are two optimal solutions and making use of duality theorem of linear programing, we can prove that these two solutions are the same. Part (b) can be proved by contradiction. 

\subsection{Insensitivity of the Solution to the Primal}


\begin{theorem}
\label{thm_insensitivity} Suppose
${\boldsymbol{x}^*}=%
\begin{pmatrix}
{\boldsymbol{x}}^{*+}\\
{\boldsymbol{x}}^{*-}%
\end{pmatrix}
$ is the optimal solution to the problem
\begin{align*}
\min~~  & f_{1}(\boldsymbol{x}^{+})\\
\text{s.t.}~~ &  (A^{\prime}- I) \boldsymbol{x}^{+} + I \boldsymbol{x}^{-} =
\boldsymbol{s}-\boldsymbol{D}\\ &  \boldsymbol{x}^{+} \ge\boldsymbol{0}, \boldsymbol{x}^{-} \ge\boldsymbol{0},
\end{align*}
where $f_{1}(\boldsymbol{x}^{+})$ is differentiable and increasing with respect to $\boldsymbol{x}^+$. Let
$f_{2}(\boldsymbol{x}^{+})$ be another differentiable and increasing function.
Then $\boldsymbol{x}^*$ is also the optimal solution to the problem
\begin{align*}
\min~~  & f_{2}(\boldsymbol{x}^{+})\\
\text{s.t.}~~ &  (A^{\prime}- I) \boldsymbol{x}^{+} + I \boldsymbol{x}^{-} =
\boldsymbol{s}-\boldsymbol{D}\\ &  \boldsymbol{x}^{+} \ge\boldsymbol{0}, \boldsymbol{x}^{-} \ge\boldsymbol{0}.
\end{align*}
\end{theorem}

To prove it, we construct the solution of the dual problem and make use of Karush-Kuhn-Tucker (KKT) conditions. See \cite{bt} for more information about KKT conditions. A detailed proof appears in Appendix \ref{appendix: thm_insensitivity}.

Although Theorem \ref{thm_insensitivity} establishes the
insensitivity of the optimal solution to a large class of nonlinear objective
functions, for the rest of the paper, our discussion is based on the primal
problem (\ref{primal}) and the dual problem (\ref{dual}) with linear objective functions.

\section{Asymptotic Behavior}

\label{sec:asym}

Now we discuss the asymptotic behavior of the failure
probability of this distribution network, which will be useful when we develop
efficient simulation algorithms for estimating the failure probability in the
next section. We will now assume fixed number $d$ of vertices in the network. We next specify the vertices' supplies and the distribution for the demands.

\review{
Let $t_i, i = 1, 2, \dots, d$, represent these $d$ locations in this network, and $T = \{t_1, t_2, \dots, t_d\}$. Suppose we have positive functions $\gamma(t),{\mu
}(t),\sigma(t)$ on $T$, and $\sigma^{2}(t,u)$ on $T\times T$. For each node $i$ with location $t_{i} \in T$,
there is a deterministic supply $s_{n}(t_{i}) \triangleq s(t_{i}) = n^{\beta
}\gamma(t_{i})$, where $\beta> 0$, $n$ is a rarity parameter, and a random
demand $D(t_{i}) \sim N({\mu}(t_{i}),\sigma^{2}(t_{i}))$, where the covariance
between the demands at two vertices with locations $t_{i}$ and $t_{j}$ is
$cov[D(t_{i}),D(t_{j})]=\sigma^{2}(t_{i},t_{j})$. }Also note that only the supply function $s(t)$ involves $n$,
not the demand function. Let $\Sigma$ be the covariance matrix of $(D(t_1), D(t_2), \ldots, D(t_d))$, which we require to be symmetric positive definite.

\review{We first introduce the little $o$ notion, which is used in the theorem that will be discussed momentarily. 

\begin{definition}
Let $f$ and $g$ be two functions defined one some subset of the real numbers. Then $f(x) = o(g(x))$ if for every $C > 0$, there exists a real number $N$ such that for all $x > N$, we have $|f(x)| < C|g(x)|$.
\end{definition}
}

We now establish a theorem that describes the asymptotic behavior of this network. More
specifically, it tells what is the most likely way in which this network
fails. This result is crucial in designing an efficient importance-sampling algorithm.




\begin{theorem}
\label{same_order} Let $L_{n}(\boldsymbol{D})$ denote the optimal value of the
dual (\ref{dual}), when the demand vector is $\boldsymbol{D}$ and the rarity
parameter is $n$. \review{Then for all $k = k_n \ge 0$ with $k_n = o(n^{\beta})$, 
\begin{align}
\lim_{n\rightarrow\infty} n^{-2\beta} \log P\{L_{n}%
(\boldsymbol{D}) > k_n\}  & = \lim_{n\rightarrow\infty}n^{-2\beta} \log
P\{\max_{i=1,\dots,d} D(t_{i})-s_n(t_{i}) > k_n\}\label{result1} \\
& = -\frac{\gamma^{2}(t^{*})}{2\sigma^{2}(t^{*})} \label{lemma},
\end{align}
where $t^{*} =\mathop{\arg\min}\limits_{t \in T} \frac{\gamma(t)}{\sigma(t)}$.}
\end{theorem}

To prove this result, we derive upper and lower bounds with the same limit $-\frac{\gamma^{2}(t^{*})}{2\sigma^{2}(t^{*})} $. The details appear in Appendix \ref{appendix: same_order}.

\section{Efficient Algorithms: Importance Sampling and Conditional Monte Carlo}

\label{sec:sim}

\subsection{Asymptotic Optimality}

Suppose $t_{i}, i=1,2,\dots, d$ are locations of $d$ vertices. When $n$ is large, the failure of this network is a rare event. To
estimate this failure probability, we develop two efficient simulation
algorithms: one based on importance sampling (IS) and the other using
conditional Monte Carlo (CMC). To evaluate the efficiency of these two
algorithms, we need to introduce a definition.

\begin{definition}
A collection $(Z_{n}: n \ge0)$ of estimators for $\rho(n)$ is said to be
asymptotically optimal if $E[Z_{n}] = \rho(n)$ and if
\[
\sup_{n>0} \frac{E(Z_{n}^{2})}{\rho(n)^{2-\epsilon}} < \infty, \forall
\epsilon> 0.
\]

\end{definition}

Asymptotic optimality also amounts to showing that
\[
\frac{\log E(Z_{n}^{2})}{2\log(\rho(n))}\rightarrow1, ~~n \rightarrow\infty.
\]

\subsection{Importance Sampling} \label{sec:IS}

We now develop an IS estimator making use of a new probability measure $Q$:
\begin{equation}
\label{Q_measure}Q\{\boldsymbol{D} \in B\} = \sum_{i=1}^{d} p(i)
P\{\boldsymbol{D} \in B | D(t_{i}) - s_n(t_{i}) > 0\},
\end{equation}
where $B \subset\mathbb{R}^{d}$ is a Borel set, and
\[
p(i) = \frac{P\{D(t_{i}) - s_n(t_{i}) > 0\}}{\sum_{j=1}^{d}P\{D(t_{j}) -
s_n(t_{j}) > 0\}}.
\]
Note that $Q$ is a mixture of $d$ measures, where the $i$-th measure in the
mixture is the conditional distribution given that the $i$-th node's demand
exceeds its supply. In other words, we force the demand to be larger than the supply for at least one node such that this network fails more often under the new measure. Since
\begin{align*}
Q\{\boldsymbol{D} \in B\}   = \frac{1}{\sum_{j=1}^{d}P\{D(t_{j}) - s_n(t_{j}) > 0\}}\sum_{i=1}^{d}
P\{\boldsymbol{D} \in B, D(t_{i}) - s_n(t_{i}) > 0\},
\end{align*}
it is easy to see that
\[
\frac{d P }{d Q } = \frac{\sum_{j=1}^{d}P\{D(t_{j}) - s_n(t_{j}) > 0\}}%
{\sum_{j=1}^{d} I\{D(t_{j}) - s_n(t_{j}) > 0 \} }.
\]

\subsubsection{Asymptotic Optimality}

We next establish the asymptotic optimality of the IS approach based on $Q$.

\begin{theorem} \label{optimal_IS}

\[
Z_{n}(\boldsymbol{D})= \frac{d P }{d Q }I\{L_{n}(\boldsymbol{D}) > k_n\} =
\frac{\sum_{j=1}^{d}P\{D(t_{j}) - s_n(t_{j}) > 0\}}{\sum_{j=1}^{d} I\{D(t_{j}) -
s_n(t_{j}) > 0 \} }I\{L_{n}(\boldsymbol{D}) > k_n\}
\]
is an asymptotically optimal estimator for $\alpha_{n}(k_n) \triangleq
P\{L_{n}(\boldsymbol{D}) > k_n\}$, \review{where $k_n = o(n^{\beta})$}.
\end{theorem}

To prove this result, we find an upper bound of $\frac{ \log E_{{Q}} [Z_n^2(\boldsymbol{D})]}{\log {P} \{L_n(\boldsymbol{D}) > k_n\} }$ with limit 2, and make use of Theorem \ref{same_order}. The proof appears in Appendix \ref{appendix: optimal_IS}.

\subsubsection{Algorithm Implementation}
\label{sec:IS_algo}

We now explain how to implement the IS algorithm.

\begin{enumerate}

\item Set $i = 1$ and let $N$ be the total number of replications to simulate.

\item Generate demand vector $\boldsymbol{D}^{(i)}$ from distribution $Q$ as
in (\ref{Q_measure}). \review{To do this, we choose a node $i$ with probability $p(i)$, and begin by generating untruncated normal variables and reject those if the demand of node $i$ does not exceed its supply. If the acceptance rate becomes too small after some iterations with escalating sample sizes, we switch to use a Gibbs sampler algorithm described in
\cite{r} to sample truncated normal variables.}

\item Calculate $Z_{n}(\boldsymbol{D}^{(i)})= \frac{\sum_{j=1}^{d}P\{D(t_{j})
- s_n(t_{j}) > 0\}}{\sum_{j=1}^{d} I\{D(t_{j}) - s_n(t_{j}) > 0 \} }%
I\{L_{n}(\boldsymbol{D}^{(i)}) > k_n\}$.

\item If $i < N$, set $i = i+1$ and go to step 2; otherwise, go to step 5.

\item Compute $\widehat{\alpha}_{n}(k_n)= (\sum_{i=1}^{N} Z_{n}(\boldsymbol{D}%
^{(i)}))/N$ as our importance-sampling estimator of $\alpha_{n}(k_n) =
P\{L_{n}(\boldsymbol{D}) > k_n\}$, and a $100(1-\delta)\%$ confidence interval
for $\alpha_{n}(k_n)$ is $(\widehat{\alpha}_{n}(k_n) \pm\Phi^{-1}(1-\delta/2)
\widehat{S}/\sqrt{N}))$, where $\widehat{S}^{2} = \big(\sum_{i=1}^{N}
(Z_{n}(\boldsymbol{D}^{(i)}) - \widehat{\alpha}_{n}(k_n))^{2}\big)/(N-1)$, and
$\Phi(\cdot)$ is the distribution function of a standard normal.
\end{enumerate}

\subsection{Conditional Monte Carlo} \label{sec:CMC}
We first briefly introduce the Conditional Monte Carlo (CMC) approach, which is a variance-reduction technique. Suppose we are interested in estimating $\alpha$, and $U$ is an unbiased estimator. According to the conditional variance formula:
$Var(U) = E[Var(U|Y)]+Var(E[U|Y])$, we have $Var(U) \ge Var(E[U|Y])$. Therefore,  using $E[U|Y]$ as an estimator may help to reduce variance. 

Now we explain how CMC is applied to our problem to estimate $\alpha(k)$. Note that the multivariate-normal random demand has polar-coordinate representation
(see
\citealp{mfe})
\begin{equation}
\label{polar}\boldsymbol{D} = \boldsymbol{\mu} + R W\boldsymbol{\Psi},
\end{equation}
where the radius $R$ satisfies $R^{2} \sim\Gamma(d/2,1/2)$, i.e., its density function $g(x) =
x^{d/2-1}e^{-x/2}(1/2)^{d/2}/\Gamma(d/2)$, $\Gamma(\cdot)$ is the gamma
function, $WW^{T} = \Sigma$, the angle $\boldsymbol{\Psi} = {\boldsymbol{z}%
}/{\Arrowvert \boldsymbol{z} \Arrowvert}$, is uniformly distributed over the
unit sphere, $\boldsymbol{z} = (z_{1},z_{2},\dots,z_{d})^{\prime}\sim N(0,I)$, and $~\Arrowvert \boldsymbol{z} \Arrowvert = \sqrt{z_{1}^{2}+z_{2}^{2}+\dots
+z_{d}^{2}}$. In addition, the radius $R$ and angle
$\mathbf{\Psi}$ are independent.

Making use of this representation, \cite{bln} developed a conditional Monte
Carlo approach for estimating $\alpha(k_n)$, along with algorithmic details on how to implement the method. However, we did not discuss the optimality of the CMC algorithm in the conference paper. We now provide such an analysis.

\subsubsection{Asymptotic Optimality}

Recall that we defined in Section~\ref{sec:asym} the deterministic
supply of node $i$ at location $t_{i}$ as $s_{n}(t_{i}) = n^{\beta}%
\gamma(t_{i})$, where $\beta> 0$ is a constant, $n$ is the rarity parameter,
and $\gamma( \, \cdot\, )$ is a fixed positive function.

\begin{theorem}\label{optimal_CMC}
{\color{white}1} For \review{$k_n = o(n^{\beta})$}, there exist $n_{0} > 0$, $c_{3} > 0$, $s^{*} > 0$, $\eta_{1} = O(n^{\beta})$, such that when $n > n_{0}$,
\begin{equation}
\label{efficiency_upper}T_{n}(\boldsymbol{\Psi}) \triangleq P\{ L_{n}%
(\boldsymbol{D})> k_n |\boldsymbol{\Psi} \} \le{P}\{ R > n^{\beta} s^{*} +
\eta_{1} \}, ~~\forall\| \boldsymbol{\Psi} \|=1,
\end{equation}
\begin{equation}
\label{efficiency lower}P\{L_{n}(\boldsymbol{D}) > k_n\} \ge c_{3} {P} \{ R >
n^{\beta} s^{*}+O(1)\}n^{-(d-1)\beta}.
\end{equation}

Also, the conditional Monte Carlo estimator $T_{n}(\boldsymbol{\Psi})$ is
asymptotically optimal.
\end{theorem}

\review{
To prove (\ref{efficiency_upper}), since the dual problem is an LP, we only need to consider the extreme points of the feasible region. Making use of the polar-coordinate representation of the random demand, we show that $P\{ L_{n}%
(\boldsymbol{D})> k_n |\boldsymbol{\Psi} \}$ is equal to the conditional probability that the radius $R$ is larger than a function of $\boldsymbol{\Psi}$, which has minimum value $n^{\beta} s^{*} +
\eta_{1}$ when $n$ is large enough. 

To prove (\ref{efficiency lower}), we show that $P\{ L_{n}%
(\boldsymbol{D})> k_n \}$ is equal to the probability that radius $R$ is larger than a function of $\boldsymbol{\Psi}$. We then find a lower bound by considering a small ball when $n$ is large enough. 

The asymptotical optimality follows since we have found an upper found of $\log \left(E[T_n^2(\boldsymbol{\Psi})]\right)$, and a lower bound of $\log \left({P}\{L_n(\boldsymbol{D}) > k_n\} \right)$, which is less than or equal to 2 when $n$ is large enough. The complete proof appears in Appendix \ref{appendix: optimal_CMC}.
}
\section{Numerical Examples}

\label{sec:exp}

Here we use the same basis for comparing the estimators using different simulation algorithms as in \cite{bln}.
Suppose we want to estimate $\alpha= E[X]$, and $X_{1}, X_{2},\dots,X_{N}$ are
independent replications of $X$. Then $\widehat{\alpha} = (\sum_{i=1}^{N}
X_{i})/N$ is an unbiased estimator of $\alpha$, and $S^{2} = (\sum_{i=1}%
^{N}(X_{i} - \widehat{\alpha})^{2})/(N-1)$ is an unbiased estimator of
$Var[X]=\sigma^{2}$, which we assume is finite. We then define the \textit{$RSE$ (relative standard error)}
as ${S}/({\sqrt{N}\widehat{\alpha}})$. To consider both the accuracy and computational efficiency when comparing different unbiased estimators, as suggested in \cite{gw}, we use the relative
measure $RSE^{2} \times CT$(Computing Time) as the criterion.

In our experiments we apply naive simulation, importance sampling, and
conditional Monte Carlo methods to different networks, and compare $RSE^{2}
\times CT$. For each example, assume $d$ locations $t_1, t_2, \cdots, t_d$ have been chosen, we give incidence matrix $H$, supply parameter $\boldsymbol{\gamma} = (\gamma(t_1), \gamma(t_2),\dots, \gamma(t_d))'$, and demand parameters $\boldsymbol{\mu}, \Sigma$. We have proven the asymptotic optimality of the IS and CMC estimators when the threshold $k$ is a constant or increases with the rarity parameter $n$. \review{Examples 1 and 2 show how failure probability changes with $n$ for constant $k_n$. Example 3 shows how failure probability changes when $k_n$ is a function of $n$, with $k = k_n = 20 \times n^{0.5}$ and $\beta = 1$}. We set the sample size $N = 10^{5}$ for all of the three examples.

\review{
We choose parameters based on the following considerations:
\begin{itemize}
\item Network size $d$: we did three experiments with networks of three different sizes $d = 3$, $10$, and $30$. We believe that a network with 30 nodes represents a sufficiently large example for actual applications. In addition, these experiments are used to compare the relative efficiency among different simulation algorithms. While larger networks take more time to simulate, we expect that the results across the methods would be similar. 

\item Incidence matrix $H$: it was chosen so that the network is irreducible.

\item Supply and demand related parameters $
\boldsymbol{\gamma}, \boldsymbol{\mu}, \Sigma$: it is not easy to obtain this information from real-life examples, so we constructed them so that failure rarely happens. 

\item Scale parameters $\beta$, rarity parameter $n$ and threshold $k$: they were chosen so that failure probability $\alpha(k_n)$ exhibits different orders of magnitude. Although our results establish asymptotic optimality of the IS and CMC estimators, the experiments consider a range of parameter values to study when $\alpha(k_n)$ is not too small so we can assess the performance.
\end{itemize}
}

\subsection{Example 1: $d=3$, fixed $k_n$}

The first example is a 3-dimensional network with the following parameters:
\begin{gather*}
H =
\begin{pmatrix}
0 & 1 & 0\\
1 & 0 & 1\\
0 & 1 & 0
\end{pmatrix}
, \quad\boldsymbol{\gamma} =
\begin{pmatrix}
3\\
1\\
13
\end{pmatrix}
, \quad\boldsymbol{\mu} =
\begin{pmatrix}
1\\
1\\
2
\end{pmatrix}
, \quad\Sigma=
\begin{pmatrix}
1 & 0.5 & 0.1\\
0.5 & 1 & 0.5\\
0.1 & 0.5 & 1
\end{pmatrix}
, \quad \beta = 1, \quad k_n = 1.
\end{gather*}

\begin{table}[h!]
\caption{{Results of Naive Simulation, IS, and CMC for $d = 3$, fixed $k_n$}. }%
\centering
\vspace{0.5cm}
\begin{tabular}
[c]{c|c|c|c|c|c|c}\hline
& \multicolumn{2}{|c}{Naive Simulation} & \multicolumn{2}{|c}{Importance
Sampling} & \multicolumn{2}{|c}{Conditional MC}\\\hline
$n$ & $\alpha(k_n)$ & $RSE^{2} \times CT$ & $\alpha(k_n)$ & $RSE^{2} \times CT$ &
$\alpha(k_n)$ & $RSE^{2} \times CT$\\\hline
1.5 & 6.77$\times10^{-2}$ & 5.04$\times10^{-2}$ & 6.76$\times10^{-2}$ &
1.59$\times10^{-2}$ & 6.69$\times10^{-2}$ & 4.35$\times10^{-2}$\\\hline
2.5 & 6.44$\times10^{-3}$ & 5.34$\times10^{-1}$ & 6.19$\times10^{-3}$ &
4.40$\times10^{-2}$ & 6.21$\times10^{-3}$ & 7.74$\times10^{-2}$\\\hline
3.2 & 6.10$\times10^{-4}$ & 5.63$\times10^{0}$ & 6.92$\times10^{-4}$ &
8.82$\times10^{-2}$ & 6.88$\times10^{-4}$ & 1.14$\times10^{-1}$\\\hline
3.9 & 8.00$\times10^{-5}$ & 4.27$\times10^{1}$ & 4.82$\times10^{-5}$ &
4.68$\times10^{-1}$ & 4.83$\times10^{-5}$ & 1.43$\times10^{-1}$\\\hline
4.5 & 0 & NaN & 3.39$\times10^{-6}$ & 1.62$\times10^{0}$ & 3.30$\times10^{-6}$
& 1.84$\times10^{-1}$\\\hline
4.9 & 0 & NaN & 4.80$\times10^{-7}$ & 7.08$\times10^{0}$ & 4.89$\times10^{-7}$
& 2.03$\times10^{-1}$\\\hline
\end{tabular}
\end{table}

\vspace{8cm}
\subsection{Example 2: $d=10$, fixed $k_n$}

The second example is a 10-dimensional network with the following parameters:

$H(i,j) = 1$ for $(i,j) = (1,2)$, $(1,3)$, $(2,1)$, $(3,4)$, $(3,8)$, $(4,5)$,$(4,7)$, $(5,6)$, $(6,7)$, $(7,8)$, $(8,9)$, $(9,10)$, $(10,1)$. All other elements of $H$ are equal to $0$.

$\boldsymbol{\gamma} = (3,5,3,3,3,3,3,3,3,15)'$, $\quad \boldsymbol{\mu} = (1,5,1,1,1,1,1,1,1,1)' $, $\quad \beta = 1,  \quad k_n = 2$.

\begin{gather*}
\Sigma=
\begin{pmatrix}
0.5 & 0.3 & 0.3 & 0.25 & 0.2 & 0.15 & 0.2 & 0.25 & 0.2 & 0.15\\
0.3 & 0.5 & 0.25 & 0.2 & 0.15 & 0.1 & 0.15 & 0.2 & 0.15 & 0.1\\
0.3 & 0.25 & 0.5 & 0.3 & 0.25 & 0.2 & 0.25 & 0.3 & 0.25 & 0.2\\
0.25 & 0.2 & 0.3 & 0.5 & 0.3 & 0.25 & 0.3 & 0.25 & 0.2 & 0.15\\
0.2 & 0.15 & 0.25 & 0.3 & 0.5 & 0.3 & 0.25 & 0.2 & 0.15 & 0.1\\
0.15 & 0.1 & 0.2 & 0.25 & 0.3 & 0.5 & 0.3 & 0.25 & 0.2 & 0.15\\
0.2 & 0.15 & 0.25 & 0.3 & 0.25 & 0.3 & 0.5 & 0.3 & 0.25 & 0.2\\
0.25 & 0.2 & 0.3 & 0.25 & 0.2 & 0.25 & 0.3 & 0.5 & 0.3 & 0.25\\
0.2 & 0.15 & 0.25 & 0.2 & 0.15 & 0.2 & 0.25 & 0.3 & 0.5 & 0.3\\
0.15 & 0.1 & 0.2 & 0.15 & 0.1 & 0.15 & 0.2 & 0.25 & 0.3 & 0.5
\end{pmatrix}
.
\end{gather*}

\begin{table}[!htb]
\caption{{Results of Naive Simulation, IS, and CMC for $d = 10$, fixed $k_n$}. }%
\centering
\vspace{0.5cm}
\begin{tabular}
[c]{c|c|c|c|c|c|c}\hline
& \multicolumn{2}{|c}{Naive Simulation} & \multicolumn{2}{|c}{Importance
Sampling} & \multicolumn{2}{|c}{Conditional MC}\\\hline
$n$ & $\alpha(k_n)$ & $RSE^{2} \times CT$ & $\alpha(k_n)$ & $RSE^{2} \times CT$ &
$\alpha(k_n)$ & $RSE^{2} \times CT$\\\hline
1.0 & 3.64$\times10^{-2}$ & 1.21$\times10^{-1}$ & 3.67$\times10^{-2}$ &
9.57$\times10^{-2}$ & 3.66$\times10^{-2}$ & 2.00$\times10^{-1}$\\\hline
1.3 & 3.05$\times10^{-3}$ & 1.39$\times10^{0}$ & 3.38$\times10^{-3}$ &
2.09$\times10^{-1}$ & 3.38$\times10^{-3}$ & 6.85$\times10^{-1}$\\\hline
1.5 & 2.10$\times10^{-4}$ & 2.00$\times10^{1}$ & 2.70$\times10^{-4}$ &
6.14$\times10^{-1}$ & 2.73$\times10^{-4}$ & 2.28$\times10^{0}$\\\hline
1.6 & 4.00$\times10^{-5}$ & 1.04$\times10^{2}$ & 3.20$\times10^{-5}$ &
2.19$\times10^{0}$ & 3.23$\times10^{-5}$ & 3.79 $\times10^{0}$\\\hline
1.7 & 0 & NaN & 4.13$\times10^{-6}$ & 1.09$\times10^{1}$ & 4.02$\times10^{-6}$
& 6.07$\times10^{0}$\\\hline
1.8 & 0 & NaN & 7.34$\times10^{-7}$ & 5.24$\times10^{1}$ & 7.26$\times
10^{-7}$ & 6.87$\times10^{0}$\\\hline
\end{tabular}
\end{table}

\newpage
\review{
\subsection{Example 3: $d=30$, $k_n$ changes with $n$}

The third example is a 30-dimensional network with the following parameters:

$H(i,i+1) = 1, i = 1,2,\dots, 29$. \quad$H(30,1) = 1$. All other elements of $H$ are equal to 0.

$\gamma(t_i) = 2, \mu(t_i) = 1, i = 1,2, \dots, 30$.  \quad $\beta = 1$,  \quad $k_n = 20 \times n^{0.5}$.

$\Sigma(i,i) = \sigma^2(t_i,t_i) = 1, i = 1,2,\dots,30$. All other elements of $\Sigma$ are equal to 0.4.

\begin{table}[tbh]
	\caption{{Results of Naive Simulation, IS, and CMC for $d = 30$, $k_n$ increases with $n$}. }%
	\centering
	\vspace{0.5cm}
	\begin{tabular}
		[c]{c|c|c|c|c|c|c}\hline
		& \multicolumn{2}{|c}{Naive Simulation} & \multicolumn{2}{|c}{Importance
			Sampling} & \multicolumn{2}{|c}{Conditional MC}\\\hline
		$n$ & $\alpha(k_n)$ & $RSE^{2} \times CT$ & $\alpha(k_n)$ & $RSE^{2} \times CT$ &
		$\alpha(k_n)$ & $RSE^{2} \times CT$\\\hline
		1.20 & 3.29$\times10^{-2}$ & 2.09$\times10^{-1}$ & 3.22$\times10^{-2}$ &
		2.94$\times10^{-1}$ & 3.23$\times10^{-2}$ & 5.96$\times10^{-1}$\\\hline
		1.50 & 2.72$\times10^{-3}$ & 2.16$\times10^{0}$ & 2.58$\times10^{-3}$ &
		1.06$\times10^{0}$ & 2.61$\times10^{-3}$ & 2.96$\times10^{0}$\\\hline
		1.70 & 2.80$\times10^{-4}$ & 2.03$\times10^{1}$ & 3.03$\times10^{-4}$ &
		3.33$\times10^{0}$ & 3.03$\times10^{-4}$ & 1.20$\times10^{1}$\\\hline
		1.95 & 1.00$\times10^{-5}$ & 5.78$\times10^{2}$ & 1.18$\times10^{-5}$ &
		2.34$\times10^{1}$ & 1.17$\times10^{-5}$ & 4.47 $\times10^{1}$\\\hline
		2.05 & 0 & NaN & 2.92$\times10^{-6}$ & 6.39$\times10^{1}$ & 3.02$\times10^{-6}$
		& 9.92$\times10^{1}$\\\hline
		2.16 & 0 & NaN & 3.83$\times10^{-7}$ & 3.07$\times10^{2}$ & 3.84$\times
		10^{-7}$ & 2.15$\times10^{2}$\\\hline
	\end{tabular}
\end{table}
}
%
%

\subsection{Discussion of Results and Comparisons Between Algorithms}

\begin{enumerate}

\item When $n$ increases, the performance of both the naive simulation and IS deteriorates quickly in terms of $RSE^{2} \times CT$.
Because we fix the number of simulations $N$, as in Example
1, 2, and 3, when $k_n$ is very large, we do not get even one observation of the
event $\{L_{n}(\boldsymbol{D})\ge k_n \}$. However, although the performance of CMC becomes worse as well, it does not deteriorate as quickly as the
other two. No matter how large $k_n$ is, we can obtain a non-zero estimate of
$\alpha(k_n)$.

\item Although both IS and CMC are asymptotically optimal, when $n$ is
small, IS performs better than CMC, as we now explain. The IS method only
needs to solve a single optimization problem to determine $Z_{n}%
(\boldsymbol{D})$ (see Section~\ref{sec:IS_algo}) in each replication $i$. In
contrast, our conditional Monte Carlo method needs to solve several
optimization problems to find the roots $R_{i}^{*}$ which equate the optimal value of the primal and the threshold $k_n$ for a fixed angle $\boldsymbol{\Psi}$ (see equation (8) in \citealp{bln}) in each
replication $i$. Thus, the added computational effort required by CMC can lead
to it performing worse than IS. However, as $n$ increases, conditional
Monte Carlo method works much better. The larger $n$ is, the bigger the
advantage CMC has compared to naive simulation. The advantage arises because
of the significant variance reduction obtained for large $n$ overwhelms
the additional computational effort. In conclusion, for a given network, IS performs best
 when $n$ is small, and CMC is better when $n$ is large.

\review{
\item We have established the asymptotic optimality of our methods as the rarity parameter $n \to \infty$. But as with any technique for which an asymptotic property has been proven, the performance for fixed $n$ when the asymptotics are not yet in effect may differ from that for large $n$, and may not outperform naive simulation. We explore this by varying $n$ in our experiments.

}
\end{enumerate}

\section{Final Comments}

\label{sec:conc}

We discuss a distribution network model with each node subjected to given fixed supply and Gaussian random demand. The unserved demand at a node is
distributed proportionally to its neighbors. The equilibrium point is determined by a linear program whose objective is to minimizing the sum of excess demands across all nodes in this network. We developed IS and CMC approaches to efficiently estimate the failure probability. Numerical results show that these two algorithms greatly outperform naive simulation, especially when the threshold $n$ is large.

We can make several extensions.
\begin{itemize}
\review{
\item 
Cost Structure:
We assume unit cost associated with pushing unit demand from one node to another. In other words, let $c_{i,j}$ be the cost by distributing unit demand from node $i$ to node $j$. Currently, $c_{i,j} = 1$ for all $(i,j) \in E$. We can generalize this setting by using a path dependent cost structure, which means $c_{i,j}$ can be different for different $(i,j)$. At the same time, the objective function of the primal problem (\ref{primal}) now becomes 
\begin{equation*}
\min~~  \sum_{i=1}^{d}\sum_{j=1}^d x_{i}^{+}a_{ij}c_{ij}.
\end{equation*}

Here, we claim that, all theorems in the paper are still valid for the generalized structure as long as $c_{i,j} > 0$ for all $(i,j) \in E$. To see this, Theorem \ref{thm_insensitivity} has generalized the cost structure for Theorems \ref{thm_feasibility} and \ref{thm_positive}. We can also prove Theorems \ref{same_order},  \ref{optimal_IS}, and \ref{optimal_CMC} with straightforward modifications. 

}
\item 
Elliptical Copula:
For CMC algorithm, note that the algorithm requires that the radial
component, $R$, is a positive continuous random variable and that we are able to calculate the root for the optimal value of the primal as a function of $R$ conditional on the angular part, $\Psi $. Therefore the conditional Monte
Carlo algorithm applies as long as the demand vector $D$ is an elliptical
copula. 

\review{
\item 
Growing Number of Nodes:
In this paper, all of our discussion focuses on a given
graph with a fixed number of nodes. We can also consider the asymptotic behavior of a graph when the number of nodes grows large. Similar properties and simulation algorithms can be developed by embedding the Gaussian vector of demands in a continuous Gaussian random field, so that Borell-TIS inequality (\citealp{at}, p.~50) can be applied in the proof of Theorem \ref{same_order}.

}
\end{itemize}

\section*{ACKNOWLEDGMENTS}
Support from NSF grants CMMI-1069064 and CMMI-1436700 is gratefully acknowledged by the first author.

The work of the third author has been supported in part by the NSF under Grants No.~CMMI-0926949, and CMMI-1200065. Any opinions, findings, and conclusions or recommendations
expressed in this material are those of the author and do not necessarily
reflect the views of the National Science Foundation.


\appendix
\section{Proof of Theorem \ref{thm_positive} } \label{appendix: thm_positive}

\begin{proof}{Proof:}
Suppose both $\boldsymbol{x}_1=\begin{pmatrix}\boldsymbol{x}_1^+\\\boldsymbol{x}_1^-\end{pmatrix}$ and $\boldsymbol{x}_2=\begin{pmatrix}\boldsymbol{x}_2^+\\\boldsymbol{x}_2^-\end{pmatrix}$ are optimal solutions. Let $\boldsymbol{d}^* = \boldsymbol{x}_1 - \boldsymbol{x}_2 = \begin{pmatrix}\boldsymbol{x}_1^+ - \boldsymbol{x}_2^+\\\boldsymbol{x}_1^- - \boldsymbol{x}_2^-\end{pmatrix}=\begin{pmatrix}\boldsymbol{d}^{*+}\\\boldsymbol{d}^{*-}\end{pmatrix}$, which is of dimension $2d$. We want to prove that $\boldsymbol{d}^* = 0$. Consider the following linear program:
\begin{align*}
(P)~~~~~~~~  \min ~~      &\boldsymbol{0}'\boldsymbol{d}      \\
\text{s.t.}~~ & \boldsymbol{1}'\boldsymbol{d}^+ = 0\\
& (A' - I) \boldsymbol{d}^+ + I \boldsymbol{d}^- = \boldsymbol{0} \\
& \boldsymbol{d} \ge \boldsymbol{e}_j,
\end{align*}
where $\boldsymbol{e}_j$ is a $2d$-dimensional vector with the $j$th element equal to 1 and other elements equal to 0. Equivalently, we write the LP $(P)$ as

\begin{align*}
\min ~~      &\boldsymbol{0}'\boldsymbol{d}      \\
\text{s.t.}~~ & B\boldsymbol{d} = \boldsymbol{0}  ~~~~~(\alpha)\\
& \boldsymbol{d} \ge \boldsymbol{e}_j, ~~~~~(\beta)
\end{align*}
where $B = \begin{pmatrix}\boldsymbol{1}' & \boldsymbol{0}'\\A' - I & I\end{pmatrix} $. Then we only need to prove the above LP is infeasible for all $ 1 \le j \le 2d$.
Consider the corresponding dual problem:
\begin{align*}
(D)~~~~~~~~ \max ~~      &\boldsymbol{\beta}'\boldsymbol{e}_j \\
\text{s.t.}~~ & B' \boldsymbol{\alpha} + \boldsymbol{\beta} = 0 \\
& \boldsymbol{\beta} \ge 0.
\end{align*}
Then, for all $m>0$, $\boldsymbol{\alpha} = \begin{pmatrix}-m \\-m\boldsymbol{1}\end{pmatrix}, \boldsymbol{\beta} = \begin{pmatrix} m\boldsymbol{1}\\ m\boldsymbol{1}\end{pmatrix}$ is a feasible solution to $(D)$ since $(I - A) \boldsymbol{1} = 0$. The value of the objective function is $m$. Due to the arbitrariness of $m$, we see that the optimal value of the dual is unbounded. Therefore, for all $1 \le j \le 2d$, the primal is infeasible. Hence, each element of d must be 0, which means that $\boldsymbol{x}_1 = \boldsymbol{x}_2$, proving part (a). \review{Note that the objective function of the LP can be of multiple forms since we only aim to prove the infeasibility, and different choice of the objective function only leads to different construction of $\boldsymbol{\alpha}$ and $\boldsymbol{\beta}$. }

To establish (b), suppose $(\boldsymbol{x}^+, \boldsymbol{x}^-)$ is the optimal solution of the primal (\ref{primal}). Suppose for some $1 \le k \le d$, both $x_k^+$ and $x_k^-$ are strictly positive, i.e., $x_k^+ > \delta$ and $x_k^- > \delta $ for some $\delta >0$. Let $\hat{x}_k^+=x_k^+ - \delta$, $\hat{x}_k^-=x_k^- - \delta$, and define a new vector $(\bar{\boldsymbol{x}}^+,\bar{\boldsymbol{x}}^-)$ as follows:
\begin{equation*}
\begin{cases}
\bar{x}_i^+ = \hat{x}_k^+, \bar{x}_i^- = \hat{x}_k^-, &~~\text{if} ~~ i=k; \\
\bar{x}_i^+ = x_i^+, \bar{x}_i^- = \bar{x_i}^+ - (D_i - s_i + \sum_{j:(j,i)\in E}\bar{x}_j^+ a_{ji}) , &~~\text{otherwise}.
\end{cases}
\end{equation*}
Then it is not hard to show that $\bar{\boldsymbol{x}} = \begin{pmatrix}\bar{\boldsymbol{x}}^+ \\ \bar{\boldsymbol{x}}^- \end{pmatrix}$ is a feasible solution to the problem (\ref{primal}). In addition, the value of the objective function at $\bar{\boldsymbol{x}}$ is strictly less than the value at $\boldsymbol{x}$, which conflicts with the optimality of $\boldsymbol{x}$. Therefore, at least one element in the pair $(x_k^+, x_k^-)$ is zero, $\forall 1 \le k \le d$. \halmos
\end{proof}

\section{Proof of Theorem \ref{thm_insensitivity} } \label{appendix: thm_insensitivity}

\begin{proof}{Proof:}
	Consider the problem
	\begin{align*}
		(P')~~~~~~~~  \min ~~      &f_1(\boldsymbol{x}^+)\\
		\text{s.t.}~~& (A' - I) \boldsymbol{x}^+ + I \boldsymbol{x}^- = \boldsymbol{s}-\boldsymbol{D} \notag &~~~~~(\boldsymbol{\alpha})\\
		& \boldsymbol{x}^+ \ge \boldsymbol{0} \notag &~~~~~(\boldsymbol{\mu})\\
		& \boldsymbol{x}^- \ge \boldsymbol{0},\notag &~~~~~(\boldsymbol{\lambda})
	\end{align*}
	Suppose $\boldsymbol{x}^*=\begin{pmatrix}\boldsymbol{x}^{*+}\\\boldsymbol{x}^{*-}\end{pmatrix}$ is the optimal solutions to $(P')$, and the Lagrange function is
	\begin{equation*}
		L(\boldsymbol{x}^*,\boldsymbol{\alpha},\boldsymbol{\mu},\boldsymbol{\lambda}) = f(\boldsymbol{x}^{*+}) + \boldsymbol{\alpha}'[(A'-I)\boldsymbol{x}^{*+} + I\boldsymbol{x}^{*-} - \boldsymbol{s} + \boldsymbol{D}] - \boldsymbol{\mu}'\boldsymbol{x}^{*+} - \boldsymbol{\lambda}'\boldsymbol{x}^{*-}.
	\end{equation*}
	Then $(\boldsymbol{x}^{*+}, \boldsymbol{x}^{*-})$ and $(\boldsymbol{\alpha}, \boldsymbol{\mu}, \boldsymbol{\lambda})$ satisfy the \emph{Karush-Kuhn-Tucher (KKT)} conditions when $f = f_1$, i.e.
	\begin{equation*}
		\begin{cases}
			&  \nabla_{\boldsymbol{x}^+}f + (A-I)\boldsymbol{\alpha} - \boldsymbol{\mu} = 0\\
			& \boldsymbol{\alpha} - \boldsymbol{\lambda} =0\\
			& x_i^{*+}\mu_i = 0, \forall i\\
			& x_i^{*-}\lambda_i = 0, \forall i\\
			&(A' - I) \boldsymbol{x}^{*+} + I \boldsymbol{x}^{*-} = \boldsymbol{s}-\boldsymbol{D}\\
			& \boldsymbol{x}^{*+} \ge 0, \boldsymbol{x}^{*-} \ge 0,  \boldsymbol{\mu} \ge 0, \boldsymbol{\lambda} \ge 0,
		\end{cases}
	\end{equation*}
	where $\nabla_{\boldsymbol{x}^+}f$ represents the gradient of $f$ with respect to $\boldsymbol{x}^+$. Now we would like to construct the dual solution vector $(\hat{\boldsymbol{\alpha}}, \hat{\boldsymbol{\mu}}, \hat{\boldsymbol{\lambda}})$, such that when $f = f_2$, $(\boldsymbol{x}^{*+}, \boldsymbol{x}^{*-})$ and $(\hat{\boldsymbol{\alpha}}, \hat{\boldsymbol{\mu}}, \hat{\boldsymbol{\lambda}})$ satisfy the above KKT conditions. Then we can claim that $(\boldsymbol{x}^{*+}, \boldsymbol{x}^{*-})$ is also the optimal solution when $f = f_2$.
	Define $\mathcal{H} = \{1 \le i \le d: x_i^{*+}>0\} $, and $\bar{\mathcal{H}} = \{1,2,\dots,d\} \backslash \mathcal{H}$. For each $i \in \mathcal{H}$, set $\hat{\mu}_i = 0$; and for each $i \in \bar{\mathcal{H}} $, set $\hat{\lambda}_i = 0$. Without loss of generality we assume that $\mathcal{H}=\{1,2,\dots,|\mathcal{H}| \}$. Let $\boldsymbol{\mu}_{\bar{\mathcal{H}}} = \{{\mu}_{|\mathcal{H}|+1}, {\mu}_{|\mathcal{H}|+2}, \dots, {\mu}_d\}$, $\boldsymbol{\lambda}_{\mathcal{H}} = \{\lambda_1, \lambda_2, \dots, \lambda_{|\mathcal{H}|}\}$, and $\boldsymbol{\xi} = \begin{pmatrix}\boldsymbol{\lambda}_{\mathcal{H}} \\ \boldsymbol{\mu}_{\bar{\mathcal{H}}} \end{pmatrix}$. Let $Q$ be a $d \times d$ diagonal matrix with the first $|\mathcal{H}|$ diagonal elements equal to $1$ and the remaining elements equal to 0. Considering the second KKT condition, the first KKT condition becomes
	\begin{align*}
		&\nabla_{\boldsymbol{x}^+}f + (A-I)\boldsymbol{\alpha} - \boldsymbol{\mu} = \nabla_{\boldsymbol{x}^+}f + (A-I)\boldsymbol{\lambda} - \boldsymbol{\mu} = \nabla_{\boldsymbol{x}^+}f + (A-I)Q\boldsymbol{\xi} - (I-Q)\boldsymbol{\xi} =0\\
		&\Rightarrow  [(I-Q)-(A-I)Q]\boldsymbol{\xi} = \nabla_{\boldsymbol{x}^+}f\\
		&\Rightarrow (I-AQ)\boldsymbol{\xi} = \nabla_{\boldsymbol{x}^+}f.
	\end{align*}
	Notice that the matrix $A$ is irreducible and stochastic. Also we claim that $Q$ cannot be the identity matrix with probability 1. To see this, suppose $Q$ is the identity matrix, in other words, $x_i^{*+} > 0, \forall 1 \le i \le d$. Note that the conclusion of Theorem \ref{thm_positive}(b) is still valid when the objective function is $f$, and the proof is exactly the same. Then $x_i^{*-} = 0, \forall 1 \le i \le d$. Adding all constraints in the primal problem (\ref{primal}) gives us $\sum_{i=1}^d D_i = \sum_{i=1}^d s_i$. But this equality holds with probability 0. Therefore, $(I-AQ)$ is invertible with probability 1,
	and $\boldsymbol{\xi} = (I-AQ)^{-1}\nabla_{\boldsymbol{x}^+}f$. Because $f$ is increasing in $\boldsymbol{x}^+$ and $(I - AQ)^{-1} \geq 0$, we have that $\boldsymbol{\xi} \ge 0$. It is obvious that $(\boldsymbol{x}^{*+}, \boldsymbol{x}^{*-})$ and $(\hat{\boldsymbol{\alpha}}, \hat{\boldsymbol{\mu}}, \hat{\boldsymbol{\lambda}}) = (Q\boldsymbol{\xi}, (I-Q)\boldsymbol{\xi}, Q\boldsymbol{\xi})$ satisfy the above KKT conditions when $f = f_2$. \Halmos
\end{proof}

\section{Proof of Theorem \ref{same_order} } \label{appendix: same_order}

\begin{proof}{Proof:}
	We will prove this result by establishing upper and lower bounds on $P\{L_n(\boldsymbol{D}) > k_n\}$. We start with deriving an upper bound. \review{ Note that $h(t)\triangleq \frac{{D}(t)-{\mu}(t)}{\sigma(t)}$ follows standard Gaussian distribution.} We first claim that
	\begin{equation}\label{subset}
	\{ L_n(\boldsymbol{D}) > k_n \} \subseteq \{ \max\limits_{i=1,\dots,d} D(t_i)-s_n(t_i) > 0 \}.
	\end{equation}
	To see this, if we assume $\max\limits_{i=1,\dots,d} D(t_i) -s_n(t_i) \le 0$, then $D(t_i) \le s(t_i), \forall i=1,2,\dots, d$. According to Theorem \ref{thm_feasibility}(b), the primal problem (\ref{primal}) is feasible, and it is easy to see that $x_i^+ = 0, x_i^- = s_n(t_i) - D(t_i) \ge 0, \forall i=1,2,\dots, d,$ is an optimal solution to the primal problem. In this case $L_n(\boldsymbol{D}) = 0$. Thus $\{ \max\limits_{i=1,\dots,d} D(t_i)-s_n(t_i) > 0 \}^c \subseteq \{ L_n(\boldsymbol{D}) > k_n \}^c$, where ``$c$'' represents the complement of a set, and Equation (\ref{subset}) is valid. Therefore,
\review{	
	\begin{align*}
		P\{L_n(\boldsymbol{D}) > k_n\} & \le P\{\max_{i=1,\dots,d} \frac{D(t_i)-s_n(t_i)}{\sigma(t_i)}   > 0\}\\
		& = P\{\max_{t \in T} (h(t)-\frac{s_n(t)-{\mu}(t)}{\sigma(t)} )> 0\}.
	\end{align*}
	Set $\hat{t} =\mathop{\arg\max}\limits_{t \in T}
	\frac{{\mu}(t)}{\sigma(t)}$. Note that when $n$ is large enough, $\frac{n^{\beta}\gamma(t^*)}{\sigma(t^*)}-\frac{{\mu}(\hat{t})}{\sigma(\hat{t})}>0$. Then
	\begin{align}\label{1}
		P\{L_n(\boldsymbol{D}) > k_n\} & \le P\{\max_{t \in T} h(t) > \frac{n^{\beta}\gamma(t^*)}{\sigma(t^*)}-\frac{{\mu}(\hat{t})}{\sigma(\hat{t})}\} \notag \\
		& \le \bar{C} \exp\{-\frac{1}{2}(\frac{n^{\beta}\gamma(t^*)}{\sigma(t^*)}-\frac{{\mu}(\hat{t})}{\sigma(\hat{t})})^2\},
	\end{align}
	 where $\bar{C}$ is some positive constant, and the last step makes use of the fact that if a random variable $X$ follows standard Gaussian distribution, then for any $x > 0$, $P\{X > x\} \le \frac{\exp\{-x^2/2\}}{x\sqrt{2\pi}}$.
	 This establishes the desired upper bound on $P\{L_n(\boldsymbol{D}) > k_n\}$.
}	

	To obtain a lower bound on the probability, define 	$g(t) \triangleq \frac{1}{\sqrt{2 \pi}} \frac{\sigma(t)}{s_n(t)-{\mu}(t)+k_n}\exp\{-\frac{(s_n(t)-{\mu}(t)+k_n)^2}{2\sigma^2(t)}\}$,	$t \in T$, 
	where $k_n \ge 0$ is some constant. We now claim that
	\begin{equation*}
		P\{L_n(\boldsymbol{D}) > k_n\}  \ge P\{\max_{i=1,\dots,d} D(t_i)-s_n(t_i) > k_n\}.
	\end{equation*}
	To see this, note that if $\max_{i=1,\dots,d} D(t_i)-s_n(t_i) > k$, then there exists some $1 \le i_0 \le d$ such that ${D}(t_{i_0})-s_n(t_{i_0}) > k_n$. Let $\boldsymbol{y}$ be the vector with the $i_0$-th element equal to 1 and the rest of the elements equal to 0. It is easy to see that $\boldsymbol{y}$ is a feasible solution to the dual problem (\ref{dual}) and $\boldsymbol{y}'(\boldsymbol{D}-\boldsymbol{s}) = D(t_{i_0})-s_n(t_{i_0}) > k$. Therefore, $L_n(\boldsymbol{D}) > k_n$. Then,
\review{
	\begin{align}
		P\{L_n(\boldsymbol{D}) > k_n\} & \ge P\{\max_{i=1,\dots,d} D(t_i)-s_n(t_i) > k_n\}\label{2}\\
		& \ge P\{D(t^*)-s_n(t^*) > k_n\}\notag\\
		& \ge  \frac{1}{\sqrt{2 \pi}}\frac{\sigma(t_{i_0})}{s_n(t_{i_0})-{\mu}(t_{i_0})+k_n}\exp\{-\frac{(s_n(t_{i_0})-{\mu}(t_{i_0})+k_n)^2}{2\sigma^2(t_{i_0})}\}\notag\\
		& = g(t^*)C,\label{3}
	\end{align}
}
	where $C$ is some positive constant, \review{and the second-to-last step applied the fact that
	if a random variable $X \sim N({\bar{\mu}}, \bar{\sigma}^{2})$, where
	$\bar{\sigma}>0$, then for all $\alpha> {\bar{\mu}}$,
	\begin{equation}
	P\{X > \alpha\} \ge\frac{1}{\sqrt{2 \pi}} \frac{\bar{\sigma}}{\alpha-{\bar
			{\mu}}}\exp\{-\frac{(\alpha-{\bar{\mu}})^{2}}{2\bar{\sigma}^{2}}\},
	\end{equation}
}
	giving us the desired lower bound on $P\{L_n(\boldsymbol{D}) > k_n\}$.
	
	Therefore, (\ref{1}), (\ref{2}), and (\ref{3}) imply for $n$ sufficiently large,
\review{
	\begin{align*}
		&\frac{1}{\sqrt{2 \pi}}\frac{\sigma(t^*)}{n^{\beta}\gamma(t^*)-{\mu}(t^*)+k_n}\exp\{-\frac{(n^{\beta}\gamma(t^*)-{\mu}(t^*)+k_n)^2}{2\sigma^2(t^*)}\}C \\
		& \le P\{\max_{i=1,\dots,d} D(t_i)-s_n(t_i) > k_n\} \\
		 &\le P\{L_n(\boldsymbol{D}) > k_n\} \le \exp\{-\frac{1}{2}(\frac{n^{\beta}\gamma(t^*)}{\sigma(t^*)}-\frac{{\mu}(\hat{t})}{\sigma(\hat{t})})^2\}\bar{C}.
	\end{align*}
	Taking logarithms, we have
	\begin{align*}
		& \log[{\frac{1}{\sqrt{2 \pi}}\frac{\sigma(t^*)}{n^{\beta}\gamma(t^*)-{\mu}(t^*)+k_n}]-\frac{(n^{\beta}\gamma(t^*)-{\mu}(t^*)+k_n)^2}{2\sigma^2(t^*)}}+\log C]\\
		& \le \log P\{\max_{i=1,\dots,d} D(t_i)-s_n(t_i) > k_n\}\\
		& \le   \log P\{L_n(\boldsymbol{D}) > k_n\} \le - \frac{1}{2}(\frac{n^{\beta}\gamma(t^*)}{\sigma(t^*)}-\frac{{\mu}(\hat{t})}{\sigma(\hat{t})})^2+\log \bar{C}.
	\end{align*}
	Because
	\begin{align*}
		&\lim_{n \rightarrow \infty}\frac{1}{n^{2\beta}}\big(\log{[\frac{1}{\sqrt{2 \pi}}\frac{\sigma(t^*)}{n^{\beta}\gamma(t^*)-{\mu}(t^*)+k_n}] -\frac{(n^{\beta}\gamma(t^*)-{\mu}(t^*)+k_n)^2}{2\sigma^2(t^*)}}+\log C\big)\\
		=& \lim_{n \rightarrow \infty}-\frac{1}{n^{2\beta}} \frac{1}{2}\frac{(n^{\beta}\gamma(t^*)-{\mu}(t^*)+k_n)^2}{\sigma^2(t^*)}
		= -\frac{\gamma^2(t^*)}{2\sigma^2(t^*)},
	\end{align*}

	it follows that
	\begin{equation*}
		\lim_{n\rightarrow \infty} n^{-2\beta} \log P\{L_n(\boldsymbol{D}) > k_n \} = \lim_{n\rightarrow\infty}n^{-2\beta} \log
		P\{\max_{i=1,\dots,d} D(t_{i})-s_n(t_{i}) > k_n\} = -\frac{\gamma^2(t^*)}{2\sigma^2(t^*)},
	\end{equation*}
	thereby verifying (\ref{result1}) and (\ref{lemma}).
}
	
\end{proof}

\section{Proof of Theorem \ref{optimal_IS} } \label{appendix: optimal_IS}

\begin{proof}{Proof:}
Let $E_Q$ denote the expectation under $Q$, so by (\ref{subset}), we have
\begin{align*}
\log E_{{Q}} [Z_n^2(\boldsymbol{D})] &= \log E_{{Q}} [ (\frac{d{P} }{d {Q} } I\{L_n(\boldsymbol{D}) > k_n\} )^2]\\
&\le   \log E_{{Q}} [ \big(\frac{d{P} }{d {Q} }  I\{ \max\limits_{i=1,\dots,d} D(t_i)-s_n(t_i) > 0\} \big)^2].
\end{align*}
Since $I\{ \max\limits_{i=1,\dots,d} D(t_i)-s_n(t_i) > 0\} = 1$ implies $\sum_{j=1}^d I\{D(t_j) - s_n(t_j) > 0 \} \ge 1$, and under measure $Q$, $\sum_{j=1}^d I\{D(t_j) - s_n(t_j) > 0 \} \ge 1$,
\begin{align*}
\frac{d{P} }{d {Q} }  I\{ \max\limits_{i=1,\dots,d} D(t_i)-s_n(t_i) > 0\} &= \frac{\sum_{j=1}^d P\{D(t_j) - s_n(t_j) > 0\}}{\sum_{j=1}^d I\{D(t_j) - s_n(t_j) > 0 \} } I\{ \max\limits_{i=1,\dots,d} D(t_i)-s_n(t_i) > 0\} \\
& \le \sum_{j=1}^d P\{D(t_j) - s_n(t_j) > 0\}.
\end{align*}
Thus
\begin{align*}
\log E_{{Q}} [Z_n^2(\boldsymbol{D})]  &\le  \log  \big(\sum_{j=1}^dP\{D(t_j) - s_n(t_j) > 0\}\big)^2 = 2 \log \sum_{j=1}^dP\{D(t_j) - s_n(t_j) > 0\}.
\end{align*}
Since
\begin{equation*}
P\{\max\limits_{i=1,\dots,d} D(t_i)-{s_n}(t_i) > 0\} \le \sum_{j=1}^dP\{D(t_j) - s_n(t_j) > 0\} \le d \times P\{\max\limits_{i=1,\dots,d} D(t_i)-{s_n}(t_i) > 0\},
\end{equation*}
we have
\begin{equation*}
\lim_{n\rightarrow \infty}\frac{\log \sum_{j=1}^dP\{D(t_j) - s_n(t_j) > 0\}}{\log P\{\max\limits_{i=1,\dots,d} D(t_i)-{s_n}(t_i) > 0\}} = 1.
\end{equation*}
Therefore,
\begin{align*}
\lim_{n\rightarrow \infty}\frac{ \log E_{{Q}} [Z_n^2(\boldsymbol{D})]}{\log {P} \{L_n(\boldsymbol{D}) > k_n\} }
\le \lim_{n\rightarrow \infty} \frac{2 \log \sum_{j=1}^dP\{D(t_j) - s_n(t_j) > 0\}}{\log P\{\max\limits_{i=1,\dots,d} D(t_i)-{s_n}(t_i) > 0\}}  \frac{\log P\{\max\limits_{i=1,\dots,d} D(t_i)-{s_n}(t_i) > 0\}}{\log {P} \{L_n(\boldsymbol{D}) > k_n\}}
= 2,
\end{align*}
where the last equation follows from Theorem \ref{same_order}. \halmos
\end{proof}

\section{Proof of Theorem \ref{optimal_CMC} } \label{appendix: optimal_CMC}

\begin{proof}{Proof:}
We first prove (\ref{efficiency_upper}).
Let $\Omega = \{\boldsymbol{y}: M\boldsymbol{y} \le \boldsymbol{1},\boldsymbol{y} \ge \boldsymbol{0} \}$ denote the feasible region of the dual problem (\ref{dual}). Then $L_n(\boldsymbol{D}) = \max \boldsymbol{y}' (\boldsymbol{\mu} +  RW\boldsymbol{\Psi} -n^{\beta}\boldsymbol{\gamma}), \boldsymbol{y} \in \Omega $, where $\boldsymbol{\gamma} = (\gamma(t_1), \gamma(t_2), \dots, \gamma(t_d))'$ as defined in Section~\ref{sec:asym}. We are interested in the failure probability, which includes two cases as we noted previously in Section~\ref{sec:model}. One case is that the primal problem is infeasible, which, according to Theorem \ref{thm_feasibility}(b), occurs if and only if when $ \boldsymbol{1}' (\boldsymbol{\mu} +  RW\boldsymbol{\Psi} -n^{\beta}\boldsymbol{\gamma}) > 0$. The other case is that the primal problem is feasible but the optimal value is greater than $k_n$. Since the dual problem is an LP, for the second case, we can focus on the extreme points of the feasible region $\Omega$. Since $k_n \geq 0$, when $\boldsymbol{y} = \boldsymbol{0}$, the optimal value is 0, so we do not have a failure. Therefore, we do not need to consider the solution $\boldsymbol{0}$ when calculating the failure probability.

Suppose $\{\tilde{\boldsymbol{y}}_i: i=1,2,\dots, m\}$ are the extreme points of $\Omega$, excluding $\boldsymbol{0}$, and we have
\begin{align*}
\{L_n(\boldsymbol{D})> k_n\} &= \{\boldsymbol{1}' (\boldsymbol{\mu} +  RW\boldsymbol{\Psi} -n^{\beta}\boldsymbol{\gamma}) > 0\} \bigcup \big[\bigcup_{i=1}^m \{\tilde{\boldsymbol{y}}_i' (\boldsymbol{\mu} +  RW\boldsymbol{\Psi} -n^{\beta}\boldsymbol{\gamma}) > k_n\}\big]\\
&= \bigcup_{i=0}^m \{\tilde{\boldsymbol{y}}_i' (\boldsymbol{\mu} +  RW\boldsymbol{\Psi} -n^{\beta}\boldsymbol{\gamma}) > k_i\},
\end{align*}
where $\tilde{\boldsymbol{y}}_0 = \boldsymbol{1}$, and
\begin{equation*}
k_i = \begin{cases}  0, & i=0 ;\\
k_n, & i=1,2,\dots, m.
\end{cases}
\end{equation*}
Let $n_1 = \max\{0,\max \limits_{i=0,1,\dots,m} \frac{ \tilde{\boldsymbol{y}}_i'\boldsymbol{\mu} - k_i}{ \tilde{\boldsymbol{y}}_i'\boldsymbol{\gamma} } \}^{1/\beta}$.Then when $n > n_1$, we have $n^{\beta} \tilde{\boldsymbol{y}}_i'\boldsymbol{\gamma} -  \tilde{\boldsymbol{y}}_i'\boldsymbol{\mu}  + k_i > 0$. Recall that $R$ is a positive random variable, so
\begin{equation*}
\tilde{\boldsymbol{y}}_i' (\boldsymbol{\mu} +  RW\boldsymbol{\Psi} -n^{\beta}\boldsymbol{\gamma}) > k_i \quad \Rightarrow \quad
\begin{cases}
R > \frac{n^{\beta} \tilde{\boldsymbol{y}}_i'\boldsymbol{\gamma}-  \tilde{\boldsymbol{y}}_i'\boldsymbol{\mu}  + k_i}{ \tilde{\boldsymbol{y}}_i' W \boldsymbol{\Psi}}, & \text{if}~~ \tilde{\boldsymbol{y}}_i' W \boldsymbol{\Psi} > 0 ;\\
R \in \emptyset, & \text{if}~~ \tilde{\boldsymbol{y}}_i' W \boldsymbol{\Psi} \le 0 .\\
\end{cases}
\end{equation*}
Define
\begin{equation*}
\Gamma_0 = \{\boldsymbol{\Psi}: \|\boldsymbol{\Psi}\|=1, \max\limits_{i=0,1,\dots,m} \tilde{\boldsymbol{y}}_i' W \boldsymbol{\Psi} > 0 \},
\end{equation*}
\begin{equation*}
M_{\boldsymbol{\Psi}} = \{i = 0,1,\dots, m:  \tilde{\boldsymbol{y}}_i' W \boldsymbol{\Psi} > 0 \}.
\end{equation*}
For $\boldsymbol{\Psi} \in \Gamma_0 $, define
\begin{equation*}
H(\boldsymbol{\Psi}, n) = \min_{i \in M_{\boldsymbol{\Psi}}} \frac{n^{\beta} \tilde{\boldsymbol{y}}_i'\boldsymbol{\gamma} -  \tilde{\boldsymbol{y}}_i'\boldsymbol{\mu} + k_i}{ \tilde{\boldsymbol{y}}_i' W \boldsymbol{\Psi}},
\end{equation*}
\begin{equation*}
S(\boldsymbol{\Psi}) = \min_{i \in M_{\boldsymbol{\Psi}}} \frac{ \tilde{\boldsymbol{y}}_i'\boldsymbol{\gamma} }{ \tilde{\boldsymbol{y}}_i' W \boldsymbol{\Psi}}, \quad
i_{\boldsymbol{\Psi}} \in \arg \min _{i \in M_{\boldsymbol{\Psi}}} \frac{ \tilde{\boldsymbol{y}}_i'\boldsymbol{\gamma} }{ \tilde{\boldsymbol{y}}_i' W \boldsymbol{\Psi}},\quad
\quad \tilde{\boldsymbol{y}}_{\boldsymbol{\Psi}} = \tilde{\boldsymbol{y}}_{i_{\boldsymbol{\Psi} } }.
\end{equation*}
It is easy to see that when $n > n_1 $,
\begin{equation}\label{Ln}
{P}\{L_n(\boldsymbol{D})> k_n \} = {P}\{ R > H(\boldsymbol{\Psi}, n)\}.
\end{equation}
In the non-trivial case when $\Gamma_0 \ne \emptyset$, there exists some $\boldsymbol{\Psi}_0 \in \Gamma_0$. Let $a = \max\limits_{i=0,1,\dots,m}  \tilde{\boldsymbol{y}}_i' W \boldsymbol{\Psi}_0 > 0$.
Define
\begin{equation*}
\Gamma_a = \{\boldsymbol{\Psi}: \|\boldsymbol{\Psi}\|=1, \max\limits_{i=0,1,\dots,m}  \tilde{\boldsymbol{y}}_i' W \boldsymbol{\Psi} \ge a \}.
\end{equation*}
Let us consider inequality (\ref{efficiency_upper}) first. We have
\begin{equation*}
T_n(\boldsymbol{\Psi}) = {P}\{ R > H(\boldsymbol{\Psi}, n)| \boldsymbol{\Psi} \} \le {P}\{ R > \inf_{ \boldsymbol{\Psi} \in \Gamma_0 } H(\boldsymbol{\Psi}, n)\} = {P}\{ R > \inf_{ \boldsymbol{\Psi} \in \Gamma_a } H(\boldsymbol{\Psi}, n)\},
\end{equation*}
and
\begin{align*}
\inf_{\boldsymbol{\Psi} \in \Gamma_a}  H(\boldsymbol{\Psi}, n) & = \inf_{\boldsymbol{\Psi} \in \Gamma_a} \min_{i \in M_{\boldsymbol{\Psi}}} \frac{n^{\beta} \tilde{\boldsymbol{y}}_i'\boldsymbol{\gamma} - \tilde{\boldsymbol{y}}_i'\boldsymbol{\mu}  + k_i}{ \tilde{\boldsymbol{y}}_i' W \boldsymbol{\Psi}}\\
& \ge  \inf_{\boldsymbol{\Psi} \in \Gamma_a} \min_{i \in M_{\boldsymbol{\Psi}}} \frac{n^{\beta} \tilde{\boldsymbol{y}}_i'\boldsymbol{\gamma} }{ \tilde{\boldsymbol{y}}_i' W \boldsymbol{\Psi}} + \inf_{\boldsymbol{\Psi} \in \Gamma_a} \min_{i \in M_{\boldsymbol{\Psi}}} \frac{-  \tilde{\boldsymbol{y}}_i'\boldsymbol{\mu}  + k_i}{ \tilde{\boldsymbol{y}}_i' W \boldsymbol{\Psi}}\\
& = n^{\beta} \inf_{\boldsymbol{\Psi} \in \Gamma_a} S(\boldsymbol{\Psi}) + \inf_{\boldsymbol{\Psi} \in \Gamma_a} \min_{i \in M_{\boldsymbol{\Psi}}} \frac{-  \tilde{\boldsymbol{y}}_i'\boldsymbol{\mu}  + k_i}{ \tilde{\boldsymbol{y}}_i' W \boldsymbol{\Psi}}.
\end{align*}
Note that both $S(\boldsymbol{\Psi})$ and $\min \limits_{i \in M_{\boldsymbol{\Psi}}} \frac{-  \tilde{\boldsymbol{y}}_i'\boldsymbol{\mu} + k_i}{ \tilde{\boldsymbol{y}}_i' W \boldsymbol{\Psi}}$ are continuous with respect to $\boldsymbol{\Psi}$ on the compact set $\Gamma_a$. Then there exist $\boldsymbol{\Psi}^* \in \Gamma_a$ and $\eta_1 = O(n^{\beta})$ such that
\begin{equation*}
\inf_{\boldsymbol{\Psi} \in \Gamma_a} S(\boldsymbol{\Psi})  = S(\boldsymbol{\Psi}^*)= \frac{ \tilde{\boldsymbol{y}}_{\boldsymbol{\Psi}^*}'\boldsymbol{\gamma} }{\tilde{\boldsymbol{y}}_{\boldsymbol{\Psi}^*}' W \boldsymbol{\Psi}^*},
\end{equation*}
\begin{equation}\label{inter_step2}
\inf_{\boldsymbol{\Psi} \in \Gamma_a}\min \limits_{i \in M_{\boldsymbol{\Psi}}} \frac{-  \tilde{\boldsymbol{y}}_i'\boldsymbol{\mu}  + k_i}{ \tilde{\boldsymbol{y}}_i' W \boldsymbol{\Psi}} = \eta_1.
\end{equation}
Therefore,
\begin{equation*}
\inf_{\boldsymbol{\Psi} \in \Gamma_a}  H(\boldsymbol{\Psi}, n) \ge n^{\beta} S(\boldsymbol{\Psi}^*) + \eta_1.
\end{equation*}
Then we have
\begin{equation*}
T_n(\boldsymbol{\Psi}) \le {P}\{ R >  n^{\beta} S(\boldsymbol{\Psi}^*) + \eta_1 \}.
\end{equation*}
Let $s^* \triangleq S(\boldsymbol{\Psi}^*)$, then (\ref{efficiency_upper}) is established.

Now we consider the inequality (\ref{efficiency lower}).
We claim that for any $\boldsymbol{\Psi}$ in $\Gamma_a$, there exists $n_2(\boldsymbol{\Psi}) > 0$ such that when $n > n_2(\boldsymbol{\Psi})$,
\begin{equation} \label{H}
H(\boldsymbol{\Psi}, n) = n^{\beta} S(\boldsymbol{\Psi})+ \frac{k_{\boldsymbol{\Psi}}- \tilde{\boldsymbol{y}}_{\boldsymbol{\Psi}}' \boldsymbol{\mu}}{ \tilde{\boldsymbol{y}}_{\boldsymbol{\Psi}}' W \boldsymbol{\Psi}},
\end{equation}
where $k_{\boldsymbol{\Psi}}$ is the $k_i$ corresponding to $\tilde{\boldsymbol{y}}_{\boldsymbol{\Psi}}$.
To see why this is true, observe that for any $i \in M_{\boldsymbol{\Psi}}$,
\begin{align*}
\lambda_i \triangleq  & n^{\beta} S(\boldsymbol{\Psi})+ \frac{k_{\boldsymbol{\Psi}}- \tilde{\boldsymbol{y}}_{\boldsymbol{\Psi}}' \boldsymbol{\mu}}{ \tilde{\boldsymbol{y}}_{\boldsymbol{\Psi}}' W \boldsymbol{\Psi}} -
\frac{n^{\beta} \tilde{\boldsymbol{y}}_i'\boldsymbol{\gamma} - \tilde{\boldsymbol{y}}_i'\boldsymbol{\mu}  + k_i}{ \tilde{\boldsymbol{y}}_i' W \boldsymbol{\Psi}} \\
=& n^{\beta} \big(S(\boldsymbol{\Psi})-\frac{ \tilde{\boldsymbol{y}}_i'\boldsymbol{\gamma} }{ \tilde{\boldsymbol{y}}_i' W \boldsymbol{\Psi}} ) +
(\frac{k_{\boldsymbol{\Psi}}- \tilde{\boldsymbol{y}}_{\boldsymbol{\Psi}}' \boldsymbol{\mu}}{ \tilde{\boldsymbol{y}}_{\boldsymbol{\Psi}}' W \boldsymbol{\Psi}} -
\frac{k_i- \tilde{\boldsymbol{y}}_i' \boldsymbol{\mu}}{ \tilde{\boldsymbol{y}}_i' W \boldsymbol{\Psi}}\big).
\end{align*}
We know that $S(\boldsymbol{\Psi})-\frac{ \tilde{\boldsymbol{y}}_i'\boldsymbol{\gamma} }{ \tilde{\boldsymbol{y}}_i' W \boldsymbol{\Psi}} \le 0$. Define
\begin{equation*}
\mathcal{I}_{\boldsymbol{\Psi}} = \{i \in M_{\boldsymbol{\Psi}}: S(\boldsymbol{\Psi})-\frac{ \tilde{\boldsymbol{y}}_i'\boldsymbol{\gamma} }{ \tilde{\boldsymbol{y}}_i' W \boldsymbol{\Psi}} = 0 \},~~~
\mathcal{I}_{\boldsymbol{\Psi}}^- = \{i \in M_{\boldsymbol{\Psi}}: S(\boldsymbol{\Psi})-\frac{ \tilde{\boldsymbol{y}}_i'\boldsymbol{\gamma} }{ \tilde{\boldsymbol{y}}_i' W \boldsymbol{\Psi}} < 0 \}.
\end{equation*}
Choose
\begin{equation*}
i_{\boldsymbol{\Psi}} \in \mathop{\arg\min} \limits_{i \in \mathcal{I}_{\boldsymbol{\Psi}}} \frac{k_{\boldsymbol{\Psi}}- \tilde{\boldsymbol{y}}_{\boldsymbol{\Psi}}' \boldsymbol{\mu}}{ \tilde{\boldsymbol{y}}_{\boldsymbol{\Psi}}' W \boldsymbol{\Psi}},
\end{equation*}
then $\lambda_i \le 0, \forall i \in \mathcal{I}_{\boldsymbol{\Psi}}$. For $i \in \mathcal{I}_{\boldsymbol{\Psi}}^-$, note that both $S(\boldsymbol{\Psi})-\frac{ \tilde{\boldsymbol{y}}_i'\boldsymbol{\gamma} }{ \tilde{\boldsymbol{y}}_i' W \boldsymbol{\Psi}} $ and $\frac{k_{\boldsymbol{\Psi}}- \tilde{\boldsymbol{y}}_{\boldsymbol{\Psi}}' \boldsymbol{\mu}}{ \tilde{\boldsymbol{y}}_{\boldsymbol{\Psi}}' W \boldsymbol{\Psi}} -
\frac{k_i- \tilde{\boldsymbol{y}}_i' \boldsymbol{\mu}}{ \tilde{\boldsymbol{y}}_i' W \boldsymbol{\Psi}}$ are bounded on $\Gamma_a$. Then there exist $\eta_2(\boldsymbol{\Psi}), \eta_3(\boldsymbol{\Psi}) > 0$, such that
\begin{equation*}
S(\boldsymbol{\Psi})-\frac{ \tilde{\boldsymbol{y}}_i'\boldsymbol{\gamma} }{ \tilde{\boldsymbol{y}}_i' W \boldsymbol{\Psi}} \le -\eta_2(\boldsymbol{\Psi}),
\end{equation*}
\begin{equation*}
-\eta_3(\boldsymbol{\Psi}) \le  \frac{k_{\boldsymbol{\Psi}}- \tilde{\boldsymbol{y}}_{\boldsymbol{\Psi}}' \boldsymbol{\mu}}{ \tilde{\boldsymbol{y}}_{\boldsymbol{\Psi}}' W \boldsymbol{\Psi}} -
\frac{k_i- \tilde{\boldsymbol{y}}_i' \boldsymbol{\mu}}{ \tilde{\boldsymbol{y}}_i' W \boldsymbol{\Psi}} \le \eta_3(\boldsymbol{\Psi}).
\end{equation*}
\review{Since $k_n = o(n^{\beta})$}, there exists $n_2(\boldsymbol{\Psi})> 0$, such that when $n > n_2(\boldsymbol{\Psi})$, $\lambda_i \le 0, \forall i \in \mathcal{I}_{\boldsymbol{\Psi}}^-$. Therefore, when $n > \max \{n_1, n_2(\boldsymbol{\Psi}^*)\}$, it follows that $\lambda_i \le 0, \forall i \in M_{\boldsymbol{\Psi}^*}$, so
\begin{equation}\label{Hstar}
H(\boldsymbol{\Psi}^*, n) = n^{\beta} S(\boldsymbol{\Psi}^*)+ \frac{k_{\boldsymbol{\Psi}^*}- \tilde{\boldsymbol{y}}_{\boldsymbol{\Psi}^*}' \boldsymbol{\mu}}{ \tilde{\boldsymbol{y}}_{\boldsymbol{\Psi}^*}' W \boldsymbol{\Psi}^*}.
\end{equation}

We also claim that there exist $c_1 > 0$, $c_2 \in \mathbb{R}$ , such that if $n > \max \{n_1, n_2(\boldsymbol{\Psi}^*) \}$, then $H(\boldsymbol{\Psi},n) - H(\boldsymbol{\Psi}^*,n ) \le (n^{\beta}c_1+c_2) \|\boldsymbol{\Psi} - \boldsymbol{\Psi}^*\|$ on $\Gamma_a$. To see this, for any $\delta > 0$ and $\boldsymbol{\theta} \in \Gamma_a$ , define $B(\boldsymbol{\theta}, \delta) = \{\boldsymbol{\Psi} \in \Gamma_a: \| \boldsymbol{\Psi}-\boldsymbol{\theta} \| \le \delta\}$. Note that there exists $\delta_1 > 0$, such that when $0 < \delta \le \delta_1$, and $n > \max \{n_1, n_2(\boldsymbol{\Psi}^*)\}$, for any $\boldsymbol{\Psi} \in B(\boldsymbol{\Psi}^*, \delta)$, we have that the index corresponding to $\tilde{\boldsymbol{y}}_{\boldsymbol{\Psi}^*}$ is in $M_{\boldsymbol{\Psi}}$, and
\begin{align*}
H(\boldsymbol{\Psi},n)- H(\boldsymbol{\Psi}^*,n)& = \min_{i \in M_{\boldsymbol{\Psi}}} \frac{n^{\beta} \tilde{\boldsymbol{y}}_i'\boldsymbol{\gamma}-  \tilde{\boldsymbol{y}}_i'\boldsymbol{\mu}  + k_i}{ \tilde{\boldsymbol{y}}_i' W \boldsymbol{\Psi} } - \frac{n^{\beta} \tilde{\boldsymbol{y}}_{\boldsymbol{\Psi}^*}'\boldsymbol{\gamma} -  \tilde{\boldsymbol{y}}_{\boldsymbol{\Psi}^*}'\boldsymbol{\mu}  + k_{\boldsymbol{\Psi}^*}}{ \tilde{\boldsymbol{y}}_{\boldsymbol{\Psi}^*}' W \boldsymbol{\Psi}^*} \\
& \le \frac{n^{\beta} \tilde{\boldsymbol{y}}_{\boldsymbol{\Psi}^*}'\boldsymbol{\gamma} -  \tilde{\boldsymbol{y}}_{\boldsymbol{\Psi}^*}'\boldsymbol{\mu}  + k_{\boldsymbol{\Psi}^*}}{ \tilde{\boldsymbol{y}}_{\boldsymbol{\Psi}^*}' W \boldsymbol{\Psi}} - \frac{n^{\beta} \tilde{\boldsymbol{y}}_{\boldsymbol{\Psi}^*}'\boldsymbol{\gamma} - \tilde{\boldsymbol{y}}_{\boldsymbol{\Psi}^*}'\boldsymbol{\mu}  + k_{\boldsymbol{\Psi}^*}}{ \tilde{\boldsymbol{y}}_{\boldsymbol{\Psi}^*}' W \boldsymbol{\Psi}^*} \\
& = (n^{\beta}\tilde{\boldsymbol{y}}_{\boldsymbol{\Psi}^*}'\boldsymbol{\gamma} -  \tilde{\boldsymbol{y}}_{\boldsymbol{\Psi}^*}'\boldsymbol{\mu}  + k_{\boldsymbol{\Psi}^*})\frac{ \tilde{\boldsymbol{y}}_{\boldsymbol{\Psi}^*}' W \boldsymbol{\Psi}^* -  \tilde{\boldsymbol{y}}_{\boldsymbol{\Psi}^*}' W\boldsymbol{\Psi}}{ \tilde{\boldsymbol{y}}_{\boldsymbol{\Psi}^*}' W \boldsymbol{\Psi}   \tilde{\boldsymbol{y}}_{\boldsymbol{\Psi}^*}' W\boldsymbol{\Psi}^*}\\
& = (n^{\beta} \tilde{\boldsymbol{y}}_{\boldsymbol{\Psi}^*}'\boldsymbol{\gamma} -  \tilde{\boldsymbol{y}}_{\boldsymbol{\Psi}^*}'\boldsymbol{\mu}  + k_{\boldsymbol{\Psi}^*})\frac{ W'\tilde{\boldsymbol{y}}_{\boldsymbol{\Psi}^*}' \boldsymbol{\Psi}^*- \boldsymbol{\Psi} }{ \tilde{\boldsymbol{y}}_{\boldsymbol{\Psi}^*}' W \boldsymbol{\Psi}  \tilde{\boldsymbol{y}}_{\boldsymbol{\Psi}^*}' W\boldsymbol{\Psi}^*}.
\end{align*}
Since $ \tilde{\boldsymbol{y}}_{\boldsymbol{\Psi}^*}' W \boldsymbol{\Psi}  \tilde{\boldsymbol{y}}_{\boldsymbol{\Psi}^*}' W\boldsymbol{\Psi}^*$ is continuous on $B(\boldsymbol{\Psi}^*, \delta)$, there exists $\delta_2 \ge 0$ such that when $0 < \delta \le \min\{\delta_1, \delta_2\}$, we have
\begin{equation*}
 \tilde{\boldsymbol{y}}_{\boldsymbol{\Psi}^*}' W \boldsymbol{\Psi} \tilde{\boldsymbol{y}}_{\boldsymbol{\Psi}^*}' W\boldsymbol{\Psi}^*
\ge  (\tilde{\boldsymbol{y}}_{\boldsymbol{\Psi}^*}' W \boldsymbol{\Psi}^*)^2 - c_0 > 0,
\end{equation*}
where $c_0$ is some positive constant.

Define $c_1 =  \tilde{\boldsymbol{y}}_{\boldsymbol{\Psi}^*}'\boldsymbol{\gamma}  \frac{\|W'\tilde{\boldsymbol{y}}_{\boldsymbol{\Psi}^*}\|}{( \tilde{\boldsymbol{y}}_{\boldsymbol{\Psi}^*}' W \boldsymbol{\Psi}^*)^2 - c_0} > 0$, $c_2 = (k_{\boldsymbol{\Psi}^*} - \tilde{\boldsymbol{y}}_{\boldsymbol{\Psi}^*}'\boldsymbol{\mu} ) \frac{\|W'\tilde{\boldsymbol{y}}_{\boldsymbol{\Psi}^*}\|}{( \tilde{\boldsymbol{y}}_{\boldsymbol{\Psi}^*}' W \boldsymbol{\Psi}^*)^2 - c_0}$. \review{Since $k_n = o(n^{\beta})$}, there exists $n_3(\boldsymbol{\Psi}^*) > 0$, such that when $n > \max \{n_1, n_2(\boldsymbol{\Psi}^*), n_3(\boldsymbol{\Psi}^*)\}$, we have $n^{\beta}c_1 + c_2 > 0$.
Therefore,
\begin{align*}
H(\boldsymbol{\Psi},n)- H(\boldsymbol{\Psi}^*,n)
 & \le (n^{\beta} \tilde{\boldsymbol{y}}_{\boldsymbol{\Psi}^*}'\boldsymbol{\gamma} - \tilde{\boldsymbol{y}}_{\boldsymbol{\Psi}^*}'\boldsymbol{\mu}  + k_{\boldsymbol{\Psi}^*})\frac{ W'\tilde{\boldsymbol{y}}_{\boldsymbol{\Psi}^*}' \boldsymbol{\Psi}^*- \boldsymbol{\Psi} }{ \tilde{\boldsymbol{y}}_{\boldsymbol{\Psi}^*}' W \boldsymbol{\Psi}  \tilde{\boldsymbol{y}}_{\boldsymbol{\Psi}^*}' W\boldsymbol{\Psi}^*}\\
 & \le  (n^{\beta} \tilde{\boldsymbol{y}}_{\boldsymbol{\Psi}^*}'\boldsymbol{\gamma} -  \tilde{\boldsymbol{y}}_{\boldsymbol{\Psi}^*}'\boldsymbol{\mu}  + k_{\boldsymbol{\Psi}^*})\frac{\|W'\tilde{\boldsymbol{y}}_{\boldsymbol{\Psi}^*}\| \|\boldsymbol{\Psi}^*- \boldsymbol{\Psi} \|}{  (\tilde{\boldsymbol{y}}_{\boldsymbol{\Psi}^*}' W \boldsymbol{\Psi}^*)^2 - c_0}\\
& = (n^{\beta}c_1 + c_2)\|\boldsymbol{\Psi}^*- \boldsymbol{\Psi} \|.
\end{align*}
So for any $\boldsymbol{\Psi} \in B(\boldsymbol{\Psi}^*, \delta)$,
\begin{equation}\label{P_bound}
H(\boldsymbol{\Psi}, n) \le H(\boldsymbol{\Psi}^*, n) + (n^{\beta}c_1 + c_2) \delta.
\end{equation}

Since $\boldsymbol{\Psi}$ is uniformly distributed over the unit sphere, which is a $(d-1)$-dimensional manifold, there exists some constant $c_3 > 0$ such that
\begin{equation*}
P \{\parallel \boldsymbol{\Psi} - \boldsymbol{\Psi}^* \parallel \le \delta  \} \ge c_3 \delta^{(d-1)}.
\end{equation*}
Let $\delta = n^{-\beta}$. By equations (\ref{Ln}) and (\ref{P_bound}), it follows that
\begin{align}
P\{L_n(\boldsymbol{D}) > k_n\}&=  P\{ R > H(\boldsymbol{\Psi}, n)\} \notag \\
& \ge {P} \{ R > H(\boldsymbol{\Psi}^*, n)+ (n^{\beta}c_1 + c_2) \delta, \parallel \boldsymbol{\Psi} - \boldsymbol{\Psi}^* \parallel \le \delta  \} \notag \\
& \ge c_3 {P} \{ R > H(\boldsymbol{\Psi}^*, n)+ (n^{\beta}c_1 + c_2) \delta\}\delta^{(d-1)} \notag \\
& = c_3 {P} \{ R > n^{\beta} S(\boldsymbol{\Psi}^*)+ \frac{k_{\boldsymbol{\Psi}^*}- \tilde{\boldsymbol{y}}_{\boldsymbol{\Psi}^*}' \boldsymbol{\mu}}{ \tilde{\boldsymbol{y}}_{\boldsymbol{\Psi}^*}' W \boldsymbol{\Psi}^*}+ (c_1 + c_2 n^{-\beta}) \} n^{-(d-1)\beta} \label{inter_step}\\
&=c_3 {P} \{ R > n^{\beta} S(\boldsymbol{\Psi}^*)+ O(1) \} n^{-(d-1)\beta} \notag.
\end{align}
Hence, we have proven (\ref{efficiency lower}).

We now establish the last part of the theorem.  By (\ref{efficiency_upper}) and (\ref{inter_step}), we have
\begin{align}
&\frac{\log \left(E[T_n^2(\boldsymbol{\Psi})]\right)}{\log \left({P}\{L(\boldsymbol{D}) > k_n\} \right)}  \le \frac{\log \left({P}^2\{ R >  n^{\beta} S(\boldsymbol{\Psi}^*) + \eta_1 \}\right) }{ \log \left(c_3 {P} \{ R > n^{\beta} S(\boldsymbol{\Psi}^*)+ \frac{k_{\boldsymbol{\Psi}^*}- \tilde{\boldsymbol{y}}_{\boldsymbol{\Psi}^*}' \boldsymbol{\mu}}{ \tilde{\boldsymbol{y}}_{\boldsymbol{\Psi}^*}' W \boldsymbol{\Psi}^*}+ (c_1 + c_2 n^{-\beta}) \} n^{-(d-1)\beta} \right) }  \notag \\
&=\frac{2 \log \left({P}\{ R >  n^{\beta} S(\boldsymbol{\Psi}^*) + \eta_1 \} \right) }{ \log c_3 + \log \left({P} \{ R > n^{\beta} S(\boldsymbol{\Psi}^*)+ \frac{k_{\boldsymbol{\Psi}^*}- \tilde{\boldsymbol{y}}_{\boldsymbol{\Psi}^*}' \boldsymbol{\mu}}{ \tilde{\boldsymbol{y}}_{\boldsymbol{\Psi}^*}' W \boldsymbol{\Psi}^*}+ (c_1 + c_2 n^{-\beta}) \} \right)-(d-1)\beta \log n }   \notag \\
& =2 \big(\frac{\log \left({P} \{ R > n^{\beta} S(\boldsymbol{\Psi}^*)+ \frac{k_{\boldsymbol{\Psi}^*}- \tilde{\boldsymbol{y}}_{\boldsymbol{\Psi}^*}' \boldsymbol{\mu}}{ \tilde{\boldsymbol{y}}_{\boldsymbol{\Psi}^*}' W \boldsymbol{\Psi}^*}+ (c_1 + c_2 n^{-\beta}) \} \right) }{\log \left({P}\{ R >  n^{\beta} S(\boldsymbol{\Psi}^*) + \eta_1 \} \right)}+\frac{\log c_3-(d-1)\beta \log n }{\log \left({P}\{ R >  n^{\beta} S(\boldsymbol{\Psi}^*) + \eta_1 \} \right) } \big)^{-1}. \label{cmc_limit}
\end{align}
Recall that $n^{\beta}c_1 + c_2 > 0$ when $n > \max \{n_1, n_2(\boldsymbol{\Psi}^*), n_3(\boldsymbol{\Psi}^*)\}$, so (\ref{inter_step2}) implies
\begin{equation*}
\eta_1 = \inf_{\boldsymbol{\Psi} \in \Gamma_a}\inf \limits_{i \in M_{\boldsymbol{\Psi}}} \frac{-  \tilde{\boldsymbol{y}}_i'\boldsymbol{\mu} + k_i}{ \tilde{\boldsymbol{y}}_i' W \boldsymbol{\Psi}} \le  \frac{k_{\boldsymbol{\Psi}^*}- \tilde{\boldsymbol{y}}_{\boldsymbol{\Psi}^*}' \boldsymbol{\mu}}{ \tilde{\boldsymbol{y}}_{\boldsymbol{\Psi}^*}' W \boldsymbol{\Psi}^*}+ (c_1 + c_2 n^{-\beta}).
\end{equation*}
Therefore,
\begin{equation*}
{P} \{ R > n^{\beta} S(\boldsymbol{\Psi}^*)+ \frac{k_{\boldsymbol{\Psi}^*}- \tilde{\boldsymbol{y}}_{\boldsymbol{\Psi}^*}' \boldsymbol{\mu}}{ \tilde{\boldsymbol{y}}_{\boldsymbol{\Psi}^*}' W \boldsymbol{\Psi}^*}+ (c_1 + c_2 n^{-\beta})\} \le {P} \{ R > n^{\beta} S(\boldsymbol{\Psi}^*)+ \eta_1\},
\end{equation*}
and
\begin{equation*}
\frac{\log {P} \{ R > n^{\beta} S(\boldsymbol{\Psi}^*)+ \frac{k_{\boldsymbol{\Psi}^*}- \tilde{\boldsymbol{y}}_{\boldsymbol{\Psi}^*}' \boldsymbol{\mu}}{ \tilde{\boldsymbol{y}}_{\boldsymbol{\Psi}^*}' W \boldsymbol{\Psi}^*}+ (c_1 + c_2 n^{-\beta})\} }{\log {P}\{ R >  n^{\beta} S(\boldsymbol{\Psi}^*) + \eta_1 \}} \ge \frac{\log {P} \{ R > n^{\beta} S(\boldsymbol{\Psi}^*)+ \eta_1\} }{\log {P}\{ R >  n^{\beta} S(\boldsymbol{\Psi}^*) + \eta_1 \}} =1.
\end{equation*}
When $n > n_4 = e^{\log c_3 / \beta(d-1)}$ , the second term inside the parentheses in (\ref{cmc_limit}) is non-negative. Then when $n > n_0 = \max\{n_1, n_2(\boldsymbol{\Psi}^*), n_3(\boldsymbol{\Psi}^*), n_4 \}$, it follows that (\ref{cmc_limit}) is bounded above by 2, thereby concluding the result. \halmos
\end{proof}




\begin{thebibliography}{}


\bibitem 
{abl}
Adler, R.J., J. H. Blanchet, J. Liu. 2012. Efficient monte carlo for high excursions of gaussian random fields. Annals of Applied Probability 22 1167-1214.

\bibitem 
{at}
Adler, R.J., J. E. Taylor. 2007. Random Fields and Geometry. Springer, New York.

\bibitem 
{a}
Asmussen, S., J. Blanchent, S. Juneja, L. Rojas-Nandayapa. 2011. Efficient simulation of tail probabilities of sums of correlated lognormals. Annals of Operations Research 189 5-23.

\bibitem 
{ag}
Asmussen, S., P. Glynn. 2007. Stochastic Simulation: Algorithms and Analysis. Springer, New York.

\bibitem 
{bt}
Bertsimas, D. and N. Tsitsiklis. 1997. Introduction to Linear Optimization. Athena Scientific, Massachusetts.

\bibitem 
{bbl}
Bienstock, D., J. Blanchet, J. Li. 2016. Stochastic models and control for electrical power line temperature. Energy Systems 7 1 173-192.

\bibitem 
{bhc}
Brechmann, E. C., Hendrich, K., and Czado, C. 2013. Conditional copula simulation for systemic risk stress testing. Insurance: Mathematics and Economics 53 3 722-732.

\bibitem 
{bln}
Blanchet, J., J. Li, M.K. Nakayama. 2011. A conditional monte carlo method for estimating the failure probability of a distribution network with random demands. Proceedings of the 2011 Winter Simulation Conference (WSC),  3832-3843.

\bibitem 
{d}
Dobson, I., B. A. Carreras, V. E. Lynch, D. E. Newman. 2007. Complex systems analysis of series of blackouts: cascading failure, critical points, and self-organization. Chaos 17 article 026103.

\bibitem 
{en}
Eisenberg, L., T.H. Noe. 2001. Systemic risk in financial systems. Management Science 47 226-249.

\bibitem 
{gw}
Glynn, P. W., W. Whitt. 1992. The asymptotic efficiency of simulation estimators. Operations Research 40 505-520.

\bibitem 
{ing}
Iyer, S. M., M.K.Nakayama, A. V. Gerbessiotis. 2009. A Markovian dependability model with cascading failures. IEEE Transactions on Computers 139 1238-1249.


\bibitem 
{k}
Kopparapu, C. 2002. Load Balancing Servers, Firewalls, and Caches. John Wiley \& Sons.

\bibitem 
{mfe}
McNeil, A. J., R. Frey, P. Embrechts. 2005. Quantitative Risk Management: Concepts, Techniques and Tools. Princeton University Press, New Jersey.

\bibitem 
{mn}
McNeil, A.J., Neslehova, J. 2009. Multivariate Archimedean copulas, $d$-monotone functions and $l_1$-norm symmetric distributions. The Annals of Statistics 37 5B 3059-3097.

\bibitem 
{pls}
Perninge, M., F. Lindskog, L. S\"{o}der. 2012. Importance sampling of injected powers for electric power system security analysis. IEEE Transactions on Power Systems 27 1 3-11. 

\bibitem 
{r}
Robert, C. P. 1995. Simulation of truncated normal variables. Statistics and Computing 5 121-125.

\bibitem
{wad}
Wadman, W.S., D.T. Crommelin, J.E. Frank. 2013. Applying a splitting technique to estimate electrical grid reliability. Proceedings of the 2013 Winter Simulation Conference, 577-588.

\bibitem
{wetal}
Wang, S. P., A. Chen, C.W. Liu, C.H. Chen, J. Shortle, J.Y. Wu. 2015. Efficient splitting simulation for blackout analysis. IEEE Transactions on Power Systems, 30 4 1775-1783.

\bibitem 
{w}
Watts, D. J. 2002. A simple model of global cascades on random networks. Proceedings of the National Academy of Sciences USA 99 5766-5771.



\end{thebibliography}



\end{document}